\documentclass[journal]{IEEEtran}

\usepackage{hyperref}
\usepackage{amsmath,amssymb,amsfonts,amsthm}
\usepackage{mathrsfs}
\usepackage{cleveref}
\usepackage{acro}
\usepackage{float}
\usepackage{algorithm}
\usepackage[noend]{algpseudocode}
\usepackage{graphicx}
\usepackage{subcaption}
\usepackage[squaren,thickqspace]{SIunits}
\usepackage{booktabs}
\usepackage[boldmath]{numprint}

\newcommand{\Real}{\ensuremath{\mathbb{R}}}

\newcommand\datamanifold{\mathbb{M}}

\newcommand{\RecSpace}{\ensuremath{X}}
\newcommand{\DataSpace}{\ensuremath{Y}}
\newcommand{\ParamSpace}{\ensuremath{\Theta}}
\newcommand{\ProdSpace}{\ensuremath{U}}

\newcommand{\approxInv}[1]{#1^{\dagger}}
\newcommand{\learnedInv}[1]{#1^{\dagger}}
\newcommand{\OpF}{\ensuremath{\mathcal{F}}}
\newcommand{\OpG}{\ensuremath{\mathcal{G}}}
\DeclareMathOperator{\IdentityOp}{\ensuremath{\text{Id}}}
\DeclareMathOperator{\OpA}{\ensuremath{\mathcal{A}}}

\DeclareMathOperator{\OpK}{\ensuremath{\mathcal{K}}}
\DeclareMathOperator{\ForwardOp}{\ensuremath{\mathcal{T}}}
\DeclareMathOperator{\ForwardOpPseudoInv}{\ensuremath{\approxInv{\ForwardOp}}}
\DeclareMathOperator{\ForwardOpInvLearned}{\ensuremath{\learnedInv{\ForwardOp}_{\!\vparam}}}
\DeclareMathOperator{\RadonTransform}{\ensuremath{\mathcal{P}}}

\DeclareMathOperator{\LogLikelihood}{\ensuremath{\mathcal{L}}}
\DeclareMathOperator{\AffineOp}{\ensuremath{\mathcal{W}}}
\DeclareMathOperator{\NonLinOp}{\ensuremath{\OpA}}
\DeclareMathOperator*{\argmin}{arg\,min}

\DeclareMathOperator{\RegOp}{\mathcal{S}}

\DeclareMathOperator{\loss}{\ensuremath{L}}
\DeclareMathOperator{\empiricalloss}{\ensuremath{\widehat{L}}}

\DeclareMathOperator{\grad}{\nabla\!}
\DeclareMathOperator{\ProxOp}{\ensuremath{prox}}

\newcommand{\stochastic}[1]{\mathsf{#1}}
\newcommand{\stsignal}{\stochastic{\signal}}
\newcommand{\stdata}{\stochastic{\data}}

\DeclareMathOperator{\Expect}{\mathbb{E}}

\newcommand{\ProbabilityMeasure}{\mu}
\newcommand{\signal}{\ensuremath{f}}

\newcommand{\signaltrue}{\signal_{\text{true}}}

\newcommand{\data}{\ensuremath{g}}
\newcommand{\noise}{\delta\data}
\newcommand{\signallearned}[1]{\ForwardOpInvLearned(#1)}

\newcommand{\vparam}{\ensuremath{\theta}}
\newcommand{\weight}{\ensuremath{w}}

\newcommand{\bias}{\ensuremath{b}}
\newcommand{\primal}{\ensuremath{f}}
\newcommand{\dual}{\ensuremath{h}}

\newcommand{\ReviewRemove}[1]{\ignorespaces}

\graphicspath{{figures/}}
\usepackage{paperacronyms}
\acuse{MRI,ADMM,GPU,CPU}
\crefname{equation}{}{}
\Crefname{equation}{}{}
\crefname{item}{}{}
\Crefname{item}{}{}
\hypersetup{%
	bookmarks=true,          
	pdfmenubar=true,	        
	pdffitwindow=false,     
	pdfstartview={FitH},    
	colorlinks=true,        
	linkcolor=red,          
	citecolor=red,          
	filecolor=magenta,      
	urlcolor=cyan,          
	pdfborder = {0,0,0}
}
\addunit{\pixel}{pixel}
\addunit{\pixels}{pixels}
\addunit{\voxel}{voxel}
\addunit{\decibel}{dB}
\addunit{\byte}{B}
\addunit{\hounsfield}{HU}

\usepackage{tikz}
\usetikzlibrary{decorations.pathreplacing,shapes,arrows,positioning}
\tikzset{
	vertex/.style = {
		circle,
		fill            = black,
		outer sep = 2pt,
		inner sep = 1pt,
	}
}
\usetikzlibrary{spy}
\usetikzlibrary{matrix}
\usetikzlibrary{overlay-beamer-styles}

\title{Learned Primal-dual Reconstruction}
\author{%
	\IEEEauthorblockN{%
		Jonas Adler\IEEEauthorrefmark{1}\IEEEauthorrefmark{2} 
		\and 
		Ozan \"{Oktem}\IEEEauthorrefmark{1}}
	\\	
	\IEEEauthorblockA{%
		\IEEEauthorrefmark{1}Department of Mathematics, KTH - Royal Institute of Technology \\
		\IEEEauthorrefmark{2}Elekta AB, Box 7593, SE-103 93 Stockholm, Sweden \\
		Email: \texttt{\{jonasadl, ozan\}@kth.se}}
	}

\begin{document}
	
	\maketitle
	
	\begin{abstract}
		We propose the Learned Primal-Dual algorithm for tomographic reconstruction. The algorithm accounts for a (possibly non-linear) forward operator in a deep neural network by unrolling a proximal primal-dual optimization method, but where the proximal operators have been replaced with convolutional neural networks. The algorithm is trained end-to-end, working directly from raw measured data and it does not depend on any initial reconstruction such as \ac{FBP}.
		
		We compare performance of the proposed method on low dose \acl{CT} reconstruction against \ac{FBP}, \ac{TV}, and deep learning based post-processing of \ac{FBP}. For the Shepp-Logan phantom we obtain $>\unit{6}{\decibel}$ \acs{PSNR} improvement against all compared methods. For human phantoms the corresponding improvement is \unit{6.6}{\decibel} over \ac{TV} and \unit{2.2}{\decibel} over learned post-processing along with a substantial improvement in the \acl{SSIM}. Finally, our algorithm involves only ten forward-back-projection computations, making the method feasible for time critical clinical applications.

	\end{abstract}
	
	
	\acresetall
	\acuse{ADMM,GPU,CPU}
	
	\section{Introduction}
	
	In an inverse problem, the goal is to reconstruct parameters characterizing the system under investigation from indirect observations. Such problems arise in several areas of science and engineering, like tomographic imaging where one seeks to visualize the interior structure of an object (2D/3D image) from indirect observations. Imaging technologies of this type, such as \ac{CT} and \ac{MRI} imaging, are  indispensable for contemporary medical diagnostics, intervention and monitoring. 
	
	There is by now a rich theory for \emph{model driven} tomographic image reconstruction. A key component is knowledge about the forward model, which describes the data formation process in absence of noise. In many applications, an explicit forward model can be derived starting out from the underlying physical principles that are utilized by the imaging modality. Another component accounts for the knowledge about the statistical properties of data and a priori information about the image to be recovered.
	
	A parallel line of development in signal processing has been the usage of deep learning for solving a wide range of tasks that can be cast as supervised learning. The success of these \emph{data driven} approaches are this far confined to tasks where knowledge about the forward model is not needed, or has little importance. As an example, an entirely data driven approach for tomographic image reconstruction applicable to clinical sized problems has yet to be demonstrated.
	
	\emph{A central question is whether one can combine elements of model and data driven approaches for solving ill-posed inverse problems.} In particular, is there a framework for incorporating the knowledge of a forward model when designing a neural network for reconstruction? 	
	The Learned Primal-Dual reconstruction method developed in this paper is such a framework. It applies to general inverse problems and it is best described in an abstract setting, so our starting point is to formalize the notion of an inverse problem.

	Mathematically, an inverse problem can be formulated as reconstructing (estimating) a signal $\signaltrue \in \RecSpace$ from data $\data \in \DataSpace$ where
	\begin{equation}\label{eq:InvProb}
	\data = \ForwardOp(\signaltrue) + \noise.  
	\end{equation}  
	Here, the reconstruction space $\RecSpace$ and the data space $\DataSpace$ are typically Hilbert Spaces, $\ForwardOp \colon \RecSpace \to \DataSpace$ is the forward operator that models how a  signal gives rise to data in absence of noise, and $\noise \in \DataSpace$ is a single sample of a $\DataSpace$-valued random variable that represents the noise component of data. 
	
	\subsection{Variational regularization}\label{sec:ClassicalReg}
	A common model driven approach for solving \cref{eq:InvProb} is to maximize the likelihood of the signal, or equivalently minimizing the negative data log-likelihood \cite{BeLaZa08}:
	\begin{equation}
	\label{eq:loglikelihood}
	\min_{\signal \in \RecSpace} \LogLikelihood\bigl(\ForwardOp(\signal), \data\bigr).
	\end{equation}
	
	This minimization is for typical choices of $\ForwardOp$ ill-posed, that is, a solution (if it exists) is unstable with respect to the data $\data$ in the sense that small changes to data results in large changes to a reconstruction. Hence, a maximum likelihood solution typically leads to over-fitting against data. 
	
	Variational regularization, also referred to as model based iterative reconstruction in medical image processing, avoids over-fitting by introducing a functional $\RegOp \colon \RecSpace \to \Real$ (regularization functional) that encodes a priori information about the true (unknown) $\signaltrue$ and penalizes unlikely solutions \cite{EnHaNe00,ScGrGrHaLe09}. Hence, instead of minimizing only the negative data log-likelihood as in \cref{eq:loglikelihood}, one now seeks to minimize a regularized objective functional by solving  
	\begin{equation}\label{eq:VarReg}
	\min_{\signal \in \RecSpace} \bigl[ \LogLikelihood\bigl(\ForwardOp(\signal), \data\bigr) + \lambda \RegOp(\signal) \bigr]
	\quad\text{for a fixed $\lambda \geq 0$.} 
	\end{equation}
	In the above, $\lambda$ (regularization parameter) governs the influence of the a priori knowledge encoded by the regularization functional against the need to fit data.

	\subsection{Machine learning in inverse problems}
	\label{sec:machine_learning_inverse}
	We here review results on learned iterative schemes, see \cite{GePerspectiveOn} for a wider review on machine learning for medical imaging and \cite{Adler17} for usage of machine learning for solving inverse problems in general.
	
	Machine learning is widely used for non-linear function approximation under weak assumptions and has recently emerged as the state of the art for several image processing tasks such as classification and segmentation. Applied to the inverse problem in \cref{eq:InvProb}, it can be phrased as the problem of finding a (non-linear) mapping $\ForwardOpInvLearned \colon \DataSpace \to \RecSpace$ satisfying the following \emph{pseudo-inverse} property:
	\[
	\ForwardOpInvLearned(\data) \approx \signaltrue
	\quad\text{whenever data $\data$ is related to $\signaltrue$ as in \cref{eq:InvProb}.}
	\]
	
	A key element in machine learning approaches is to parametrize the set of such pseudo-inverse operators by a parameter $\vparam \in \ParamSpace$ where $\ParamSpace$ is some parameter space and the main algorithmic complication is to select an appropriate structure of $\ForwardOpInvLearned$ such that, given appropriate training, the pseudo-inverse property is satisfied as well as possible.
	
	In the context of tomographic reconstruction, three main research directions have been proposed. The first is so called learned post-processing or learned denoisers. Here, the learned reconstruction operator is of the form
	\[
	\ForwardOpInvLearned = \Lambda_{\vparam} \circ \ForwardOpPseudoInv
	\]
	where $\Lambda_{\vparam} : \RecSpace \to \RecSpace$ is a learned post-processing operator and $\ForwardOpPseudoInv : \DataSpace \to \RecSpace$ is some approximate pseudo-inverse, e.g. given by \ac{FBP} in \ac{CT} reconstruction. This type of method is relatively easy to implement, given that the pseudo-inverse can be applied off-line, before the learning is performed, which reduces the learning to inferring an $\RecSpace \to \RecSpace$ transformation. This has been investigated by several authors \cite{2017arXiv170200288C, JiMcFrUn16, WaveNet}.
	
	Another method is to learn a regularizer and use this regularizer in a classical variational reconstruction scheme according to \cref{eq:VarReg}. Examples of this include dictionary learning \cite{QXuDictionary}, but several alternative methods have been investigated, such as learning a variational auto-encoder \cite{2017arXiv170403488M} or using a cascade of wavelet transforms (scattering transform) \cite{MallatInverse}.
	
	Finally, some authors investigate learning the full reconstruction operator, going all the way from data to reconstruction. Doing this in one step is typically very computationally expensive and does not scale to the data sizes encountered in tomographic reconstruction. Instead, \emph{learned iterative schemes} have been studied. These schemes resemble classical optimization methods used for tomographic reconstruction  but use machine learning to find the best update in each iteration given the last iterate and results of applying the forward operator and its adjoint as input.
	
	
	One of the first works on learned iterative schemes is \cite{YaSuLiXu16}, which learns an \ac{ADMM}-like scheme for \ac{MRI} reconstruction. A further development along this lines is given in \cite{PuWe17}, which learns over a broader class of schemes instead of \ac{ADMM}-type of schemes. The application is to finite dimensional inverse problems typically arising in image restoration. This approach was  in \cite{Adler17} further extended to non-linear forward operators in to the infinite dimensional setting, which also applies learned iterative schemes to (non-linear, pre-log) \ac{CT}. Similar approaches for \ac{MRI} reconstruction have also been considered \cite{PockVariational, StanfordMRI}. Here, the situation is simpler than \ac{CT} reconstruction since the forward operator is approximated by a Fourier transform, i.e. \ac{MRI} reconstruction amounts to inverting the Fourier transform.
	
	Given a structure of $\ForwardOpInvLearned$, the "learning" part refers to choosing an "optimal" set of parameters $\vparam$ given some \emph{training data}, where the concept of optimality is typically quantified through a loss functional that measures the quality of a learned pseudo-inverse $\ForwardOpInvLearned$.
	
	To define this loss functional, consider a $(\RecSpace \times \DataSpace)$--valued random variable $(\stsignal,\stdata)$ with joint probability distribution $\ProbabilityMeasure$. This could be e.g. the probability distribution of human bodies and corresponding noisy tomographic data.
	We define the optimal reconstruction operator as the one whose reconstructions have the lowest average mean squared distance to the true reconstructions, where the average is taken w.r.t. $\ProbabilityMeasure$. Finding this reconstruction operator is then given by selecting the parameters $\vparam \in \ParamSpace$ so that the loss functional $\loss(\vparam)$ is minimized:
	\begin{equation}
	\label{eq:LossStandard}
	\loss(\vparam) :=
	\Expect_{(\stsignal,\stdata) \sim \ProbabilityMeasure} 
	\Bigl[
	\bigl\Vert \signallearned{\stdata} - \stsignal \bigr\Vert_\RecSpace^2
	\Bigr].
	\end{equation}
	
	However, in practice we often do not have access to the probability distribution $\ProbabilityMeasure$ of the random variable $(\stsignal,\stdata)$. Instead, we know a finite set of samples $(\data_1, \signal_1), \dots, (\data_N, \signal_N)$. In this setting, we replace $\ProbabilityMeasure$ in \cref{eq:LossStandard} with its empirical counterpart, so the loss function is replaced with the empirical loss
	\begin{equation}
	\label{eq:LossEmpirical}
	\empiricalloss(\vparam) :=
	\frac{1}{N}
	\sum_{i=1}^N
	\Bigl[
	\bigl\Vert \signallearned{\data_i} - \signal_i \bigr\Vert_\RecSpace^2
	\Bigr].
	\end{equation}
	Since our main goal is to minimize the loss functional, our practical goal is thus two-fold: we want to find a learned reconstruction scheme that minimizes the empirical loss, and we also want it to generalize to new, unseen data.

	\section{Contribution and overview of paper}
	
	This paper proposes the \emph{Learned Primal-Dual} algorithm, a general framework for solving inverse problems that combines deep learning with model based reconstruction. The proposed learned iterative reconstruction scheme involves \acp{CNN} in both the reconstruction and data space, and these are connected by the forward operator and its adjoint. We train the networks to minimize the mean squared error of the reconstruction and demonstrate that this achieves very high performance in \ac{CT} reconstruction, surpassing recent learning based methods on both analytical and human data.
	
	We emphasize that we learn the whole reconstruction operator, mapping data to reconstruction, and not just a post-processing nor only the prior in isolation.
	
	In addition, we make all of our code and learned parameters open source so that the community can reproduce the results and apply the methods to other inverse problems \cite{SourceCode}.
	
	\section{The Learned Primal-Dual algorithm}\label{sec:GeneralAlg}
	
	We here introduce how primal-dual algorithms can be learned from data and how this can be used to solve inverse problems.
	
	\subsection{Primal-Dual optimization schemes}
	
	In imaging, the minimization in \cref{eq:VarReg} is a large scale optimization problem, which traditionally  has been addressed using gradient based methods such as gradient descent or its extensions to higher order derivatives, e.g. quasi-Newton or Newton methods. 
	However,  many regularizers of interest result in a non-differentiable objective functional, so gradient based methods are not  applicable. 	
	A common approach to handle this difficulty is to consider a smooth approximation, which however introduces additional parameters and gives non-exact solutions. 
	
	An alternative approach is to use methods from non-smooth convex optimization.
	\emph{Proximal} methods have been developed in order to directly work with non-smooth objective functionals. Here, a proximal step replaces the gradient step. The simplest example of such an algorithm is the proximal point algorithm for minimizing an objective functional $\OpG : \RecSpace \to \Real$. It can be seen as the proximal equivalent of the gradient descent scheme and is given by
	\begin{equation}
	\primal_{i + 1} = \ProxOp_{\tau \OpG}(\primal_i)
	\label{eq:proxpt}
	\end{equation}
	where $\tau \in \Real^+$ is a step size and the \emph{proximal operator} is defined by
	\begin{equation}
	\ProxOp_{\tau \OpG}(\primal) = \argmin_{\primal' \in \RecSpace} \Bigl[ \OpG(\primal') + \frac{1}{2 \tau} \bigl\Vert \primal' - \primal \bigr\Vert_\RecSpace^2 \Bigr]
	\label{eq:proxdef}
	\end{equation}
	
	While this algorithm could, in theory, be applied to solve \cref{eq:VarReg} it is rarely used directly since \cref{eq:proxdef} does not have a closed form solution.
	Proximal primal-dual schemes offer a work around. In these schemes, an auxiliary \emph{dual} variable in the range of the operator is introduced and the \emph{primal} ($\primal \in \RecSpace$) and dual variables are updated in an alternating manner.
	
	One well known primal-dual scheme is the \ac{PDHG} algorithm \cite{ChPo10}, also known as the Chambolle-Pock algorithm, with a recent extension to non-linear  operators \cite{Va14}. The scheme (\cref{alg:cp}) is adapted for minimization problems with the following structure:
	\begin{equation}
	\label{eq:cp_problem_def}
	\min_{\primal \in \RecSpace} \Bigl[ \OpF\bigl( \OpK(\primal)\bigr) + \OpG(\primal) \Bigr]
	\end{equation}
	where $\OpK : \RecSpace \to \ProdSpace$ is a (possibly non-linear) operator, $\ProdSpace$ is a Hilbert space and $\OpF \colon \ProdSpace \to \Real$ and $\OpG \colon \RecSpace \to \Real$ are functionals on the dual/primal spaces. Note that \cref{eq:VarReg} is a special case of \cref{eq:cp_problem_def} if we set $\OpF := \LogLikelihood(\cdot, g)$, $\OpK := \ForwardOp$ and $\OpG := \RegOp$.
	\begin{algorithm}[h]
		\caption{Non-linear \acl{PDHG}}\label{alg:cp}
		\begin{algorithmic}[1]
			\State Given: $\sigma, \tau > 0$ s.t. $\sigma \tau \lVert \OpK \rVert^2 < 1$, $\gamma \in [0, 1]$ and $\primal_0 \in \RecSpace$, $\dual_0 \in \ProdSpace$. 
			\For{$i = 1, \dots$}
			\State $\dual_{i + 1} \gets \ProxOp_{\sigma \OpF^*}\bigl(\dual_{i} + \sigma \OpK(\bar{\primal}_{i}) \bigr)$ 
			\State $\primal_{i + 1} \gets \ProxOp_{\tau \OpG}\bigl(\primal_{i} - \tau [ \partial \OpK(\primal_{i})]^* (\dual_{i + 1}) \bigr)$
			\State $\bar{\primal}_{i + 1} \gets \primal_{i + 1} + \gamma (\primal_{i + 1} - \primal_{i})$
			\EndFor
		\end{algorithmic}
	\end{algorithm}
	
	In \cref{alg:cp}, $\OpF^*$ denotes the Fenchel conjugate of $\OpF$, $\dual \in \ProdSpace$ is the dual variable and $[ \partial \OpK(\primal_{i})]^* : \ProdSpace \to \RecSpace$ is the adjoint of the (Fr\'echet) derivative of $\OpK$ in point $\primal_{i}$. 
	
	\paragraph*{Example: \Ac{TV} regularized \ac{CT}}
	The \ac{PDHG} method has been widely applied to \ac{CT} \cite{SidkyCT}. In \ac{CT}, the forward operator is given by the ray-transform $\RadonTransform : \RecSpace \to \DataSpace$, which integrates the signal over a set of lines $\datamanifold$ given by the acquisition geometry. Hence, elements in $\DataSpace$ are functions on lines
	\[
	\RadonTransform(\signal)(\ell) := \int_{\ell} \signal(x) dx \quad \text{for $\ell \in \datamanifold$.}
	\]
	and the adjoint of the derivative is the back-projection \cite{Ma06}.	 
	
	A typical example of variational regularization in imaging is \ac{TV} regularization, which applies to signals that are represented by scalar functions of bounded variation. The corresponding regularization functional is then given as the 1-norm of the gradient magnitude, i.e.
	$\RegOp(\signal) := \Vert \grad \signal \Vert_1$,  $\nabla : \RecSpace \to \RecSpace^d$, $d$ is the dimension of the space.

	The \ac{PDHG} method can be used to solve the \ac{TV} regularized \ac{CT} optimization problem
	\[
	\min_{\signal \in \RecSpace} \lVert \RadonTransform(\signal) - \data \rVert_2^2 + \lambda \lVert \nabla \signal \rVert_1.
	\]
	Since the proximal of $\signal \mapsto \lVert \nabla \signal \rVert_1$ is hard to compute, the following identification is better suited for recasting the above into \eqref{eq:cp_problem_def}:
	\[
	\OpK \colon \RecSpace \to \DataSpace \times \RecSpace^d
	\quad
	\text{as}
	\quad
	\OpK(\primal) := \bigl[ \RadonTransform(\primal), \nabla \primal \bigr],
	\]
	\[ 
	\OpF\bigl([\dual^{(1)}, \dual^{(2)}] \bigr) := \lVert \dual^{(1)} - \data \rVert_2^2 + \lVert \dual^{(2)} \rVert_1
	\quad\text{and}\quad
	\OpG(\primal) := 0.
	\]

	\subsection{Learned \ac{PDHG}}
	
	The aim is to derive a learned reconstruction scheme inspired by \ac{PDHG}, \cref{alg:cp}.
	We follow the observation in \cite{FlexISP, EladRED}, that proximal operators can be replaced by other operators that are not necessarily proximal operators. The aforementioned publications replace a proximal operator with a denoising operator such as Block Matching 3D (BM3D). Our idea is to replace the proximal operators by \emph{parametrized} operators where the parameters are learned from training data, resulting in a learned reconstruction operator.
	
	In order to make the learned reconstruction operator well defined and implementable on a computer we also need to select a stopping criterion. Choosing a proper stopping criterion is an active research area, but for simplicity and usability we use a fixed number of iterates. By selecting a fixed number of iterations, the computation budget is also fixed prior to training, which is a highly desirable property in time critical applications.
	
	\Cref{alg:learned_cp} below outlines the resulting variant of the \ac{PDHG} algorithm with $I$ iterations in which the primal proximal has been replaced by a learned proximal, $\Gamma_{\vparam^d}$ and the dual proximal by a learned proximal $\Lambda_{\vparam^p}$.  Note that in this article we consider only a single forward model and no regularizing operator, so we have $\OpK = \ForwardOp$, $\ProdSpace = \DataSpace$, but we give the algorithm in full generality for completeness.
	
	\begin{algorithm}[H]
		\caption{Learned \ac{PDHG}}\label{alg:learned_cp}
		\begin{algorithmic}[1]
			\State Initialize $\primal_0 \in \RecSpace, \dual_0 \in \ProdSpace$
			\For{$i = 1, \dots, I$}
			\State $\dual_{i + 1} \gets
			\Gamma_{\vparam^d}\bigl(\dual_{i} + \sigma \OpK(\bar{\primal}_{i}), \data \bigr)$ 
			\State $\primal_{i + 1} \gets
			\Lambda_{\vparam^p}\bigl(\primal_{i} - \tau [ \partial \OpK(\primal_{i})]^* (\dual_{i + 1}) \bigr)$
			\State $\bar{\primal}_{i + 1} \gets \primal_{i + 1} + \theta (\primal_{i + 1} - \primal_{i})$
			\EndFor
			\State \Return $\primal_I^{(1)}$
		\end{algorithmic}
	\end{algorithm}
	In \cref{alg:learned_cp}, there are several parameters that need to be selected. These are the parameters of the dual proximal, $\vparam^d$, the primal proximal, $\vparam^p$, the step lengths, $\sigma$, $\tau$ and the overrelaxation parameter, $\theta$. In a learned \ac{PDHG} algorithm these would all be infered, \emph{learned}, from training data.
	
	We implemented this algorithm and show its performance in the results section. While the performance was comparable to traditional methods, it did not improve upon the state of the art in deep learning based image reconstruction.
	
	\subsection{Learned Primal-Dual}
	
	To gain substantial improvements, guided by recent advances in machine learning, the following modifications to the learned \ac{PDHG} algorithm were done.
	
	\begin{itemize}
		\item Following \cite{PuWe17, Adler17}, extend the primal space to allow the algorithm some "memory" between the iterations.
		\[
		\primal = [\primal^{(1)}, \primal^{(2)}, \dots, \primal^{(N_{\text{primal}})}] \in \RecSpace^{N_{\text{primal}}}
		\]
		Similarly extend the dual space $\ProdSpace$ to $\ProdSpace^{N_{\text{dual}}}$.
		\item Instead of explicitly enforcing updates of the form $\dual_{i} + \sigma \OpK \bigl(\bar{\primal}_{i}\bigr)$, allow the network to learn how to combine the previous update with the result of the operator evaluation.
		\item Instead of hard-coding the over-relaxation $\bar{\primal}_{i + 1} \gets \primal_{i + 1} + \theta (\primal_{i + 1} - \primal_{i})$, let the network  freely learn in what point the forward operator should be evaluated.
		\item Instead of using the same learned proximal operators in each iteration allow them to differ. This increases the size of the parameter space but it also notably improves reconstruction quality.
	\end{itemize}
	
	The above modifications result in a new algorithm, henceforth called the \emph{Learned Primal-Dual} algorithm, that is outlined in \cref{alg:learned_pd}.
	\begin{algorithm}[H]
		\caption{Learned Primal-Dual}\label{alg:learned_pd}
		\begin{algorithmic}[1]
			\State Initialize $\primal_0 \in \RecSpace^{N_{\text{primal}}}, \dual_0 \in \ProdSpace^{N_{\text{dual}}}$
			\For{$i = 1, \dots, I$}
			\State $\dual_i \gets
			\Gamma_{\vparam_i^d}\bigl(\dual_{i - 1}, \OpK(\primal_{i-1}^{(2)}), \data\bigr)$
			\State $\primal_i \gets 
			\Lambda_{\vparam_i^p}\bigl(\primal_{i - 1}, [\partial\!\OpK(\primal_{i-1}^{(1)})]^*(\dual_i^{(1)}) \bigr)$
			\EndFor
			\State \Return $\primal_I^{(1)}$
		\end{algorithmic}
	\end{algorithm}
	
	\subsubsection{Choice of starting point}
	
	In theory, the Learned Primal-Dual algorithm can be used with any choice of starting points $\primal_0$ and $\dual_0$.	
	The most simple starting point, both from a conceptual and computational perspective, is zero-initialization
	\begin{align}
	\label{eq:init_zero}
	\primal_0 =&\ [0, 0, \dots, 0] \nonumber \\
	\dual_0 =&\ [0, 0, \dots, 0]
	\end{align}
	where $0$ is the zero element in the primal or dual space.
	
	In cases where a good starting guess is available, it would make sense to use it. One such option is to assume that there exists a pseudo-inverse $\ForwardOpPseudoInv$, e.g. \ac{FBP} for \ac{CT}. For the dual variable, the data $\data$ enters into each iterate so there is no need for a good initial guess. This gives the starting point
	\begin{align}
	\label{eq:init_pseudoinverse}
	\primal_0 =&\ [\ForwardOpPseudoInv(\data), \ForwardOpPseudoInv(\data), \dots, \ForwardOpPseudoInv(\data)] \nonumber \\
	\dual_0 =&\ [0, 0, \dots, 0]
	\end{align}
	
	In our tests we found that providing the Learned Primal-Dual algorithm with such an initial guess marginally decreased training time, but did not give better final results. Given that using the pseudo-inverse $\ForwardOpPseudoInv$ adds more complexity by making the learned reconstruction operator depend on an earlier reconstruction, we \emph{report values only from zero-initialization}.
	
	\subsection{Connection to variational regularization}
	
	We note that by selecting $N_{\text{primal}} = 2$ and $N_{\text{dual}} = 1$ the Learned Primal-Dual algorithm naturally reduces to the classical \ac{PDHG} algorithm by making the following choices:
	\begin{multline*}
	\Gamma_{\vparam_i^d}\bigl(\dual, \OpK(\primal^{(2)}), \data \bigr)
	=\
	\ProxOp_{\sigma \OpF_\data^*}\bigl(\dual + \sigma \OpK(\primal^{(2)}) \bigr)
	\\
	\shoveleft{\Lambda_{\vparam_i^p}\left(
		\begin{bmatrix}
		\primal^{(1)} \\
		\primal^{(2)}
		\end{bmatrix}
		,
		[\partial\!\OpK(\primal^{(1)})]^*(\dual) \right)
		=\ } \\
	\begin{bmatrix}
	\ProxOp_{\tau \OpG}\bigl(\primal^{(1)} - \tau [\partial\!\OpK(\primal^{(1)})]^*(\dual) \bigr)
	\\
	(1 + \theta) \ProxOp_{\tau \OpG}\bigl(\primal^{(1)} - \tau [\partial\!\OpK(\primal^{(1)})]^*(\dual) \bigr) - \theta \primal^{(1)}
	\end{bmatrix}
	.
	\end{multline*}
	Even if the learned proximal operators do not have explicit access to the proximals, the universal approximation property of neural networks \cite{UniversalApproximation} guarantees that given sufficient training data these equalities can be approximated arbitrarily well.
	
	A wide range of other optimization schemes can also be seen as special cases of the Learned Primal-Dual algorithm. For example, the gradient descent algorithm with step-length $\alpha$ for solving \cref{eq:cp_problem_def} is given by
	\begin{equation*}
	\signal_{i + 1} = 
	\signal_{i} - \alpha \Bigl([\partial\!\OpK(\primal_i)]^* \bigl([\nabla \OpF](\OpK(\primal_i))\bigr)
	+ [\nabla \OpG](\primal_i) \Bigr)
	\end{equation*}
	and can be obtained by selecting 
	\begin{multline*}
	\Gamma_{\vparam_i^d}\bigl(\dual, \OpK(\primal^{(2)}), \data \bigr)
	=\
	[\nabla \OpF_\data](\OpK(\primal^{(2)}))
	\\
	\shoveleft{
		\Lambda_{\vparam_i^p}\left(
		\begin{bmatrix}
		\primal^{(1)} \\
		\primal^{(2)}
		\end{bmatrix}
		,
		[\partial\!\OpK(\primal^{(1)})]^*(\dual) \right)
		=
	}
	\\
	\begin{bmatrix}
	\primal^{(1)} - \alpha \bigl([\nabla \OpG](\primal^{(1)}) + [\partial\!\OpK(\primal^{(1)})]^*(\dual) \bigr)
	\\
	\primal^{(1)} - \alpha \bigl([\nabla \OpG](\primal^{(1)}) + [\partial\!\OpK(\primal^{(1)})]^*(\dual) \bigr)
	\end{bmatrix}
	.
	\end{multline*}
	More advanced gradient based methods such as Limited memory BFGS are likewise sub-cases obtained by appropriate choices of learned proximal operators.
	
	In summary, the Learned Primal-Dual algorithm contains a wide range of optimization schemes as special cases. If the parameters are appropriately selected, then the proposed algorithm should always perform \emph{at least as well} as current variational regularization schemes given the \emph{same stopping criteria}, which here is a fixed number of iterates.

	\section{Implementation and evaluation}\label{sec:Eval}
	
	We evaluate the algorithm on two low dose \ac{CT} problems. One simplified using analytical phantoms based on ellipses and one  with a more realistic forward model and human phantoms. We briefly describe these test cases and how we implemented the Learned Primal-Dual algorithm. We also describe the methods we compare against.
	
	\subsection{Test cases}
	
	\paragraph{Ellipse phantoms}
	This problem is identical to \cite{Adler17} and we restate it briefly.
	Training data is randomly generated ellipses on a \unit{128 \times 128}{\pixel} domain. The forward operator is the ray transform and hence $\ForwardOp = \RadonTransform$.
	
	The projection geometry was a sparse 30 view parallel beam geometry with 182 detector pixels. 5\% additive Gaussian noise was added to the projections. Since the forward operator is linear, the adjoint of the derivative is simply the adjoint, which for the ray transform is the back-projection 
	\[
	[\partial \ForwardOp(\signal)]^*
	=
	\RadonTransform^*.
	\]
	\paragraph{Human phantoms}
	In order to evaluate the algorithm on a clinically realistic use-case we consider reconstruction of simulated data from human abdomen \ac{CT} scans as provided by Mayo Clinic for the AAPM Low Dose CT Grand Challengfe \cite{AAPMLowDose}. The data includes full dose \ac{CT} scans from 10 patients, of which we used 9 for training and 1 for evaluation. We used the \unit{3}{\milli\meter} slice thickness reconstructions, resulting in 2168 training images, each $\unit{512 \times 512}{\pixel}$ in size. Thus, given that we minimize the pointwise error, the total number of data-points is $512^2 \times 2168 \approx 5 \times 10^8$.

\begin{figure}[t]
		\centering	
		\begin{subfigure}[t]{.48\linewidth}	
			\includegraphics[width=\linewidth, trim={23mm 17mm 32mm 6mm}, clip]{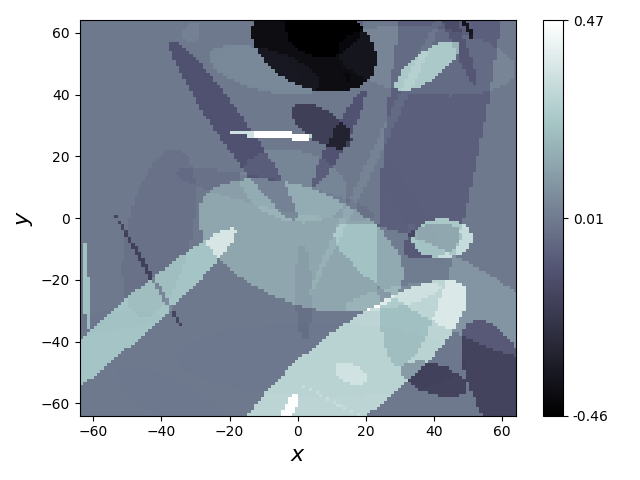}
			\caption{Ellipse training phantom}
		\end{subfigure}
		\begin{subfigure}[t]{.48\linewidth}
			\includegraphics[width=\textwidth, trim={23mm 17mm 32mm 6mm}, clip]{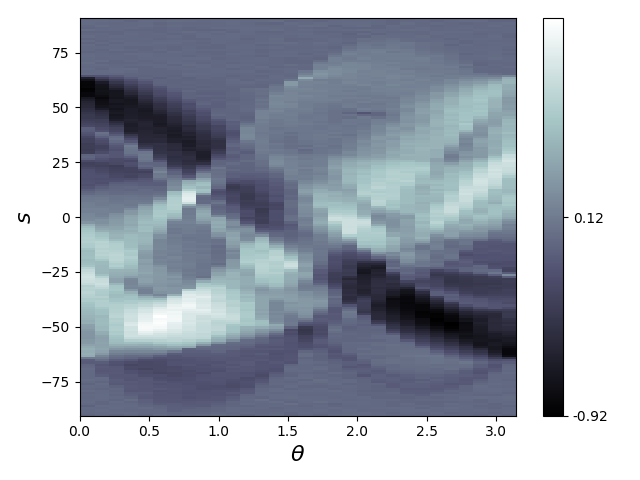}
			\caption{Ellipse training sinogram}
		\end{subfigure}
		\\
		\begin{subfigure}[t]{.48\linewidth}	
			\includegraphics[width=\linewidth, trim={23mm 17mm 32mm 6mm}, clip]{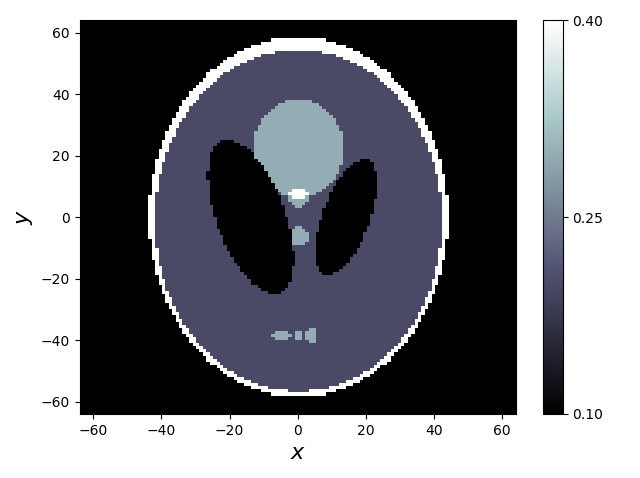}
			\caption{Shepp-Logan phantom}
		\end{subfigure}
		\begin{subfigure}[t]{.48\linewidth}
			\includegraphics[width=\textwidth, trim={23mm 17mm 32mm 6mm}, clip]{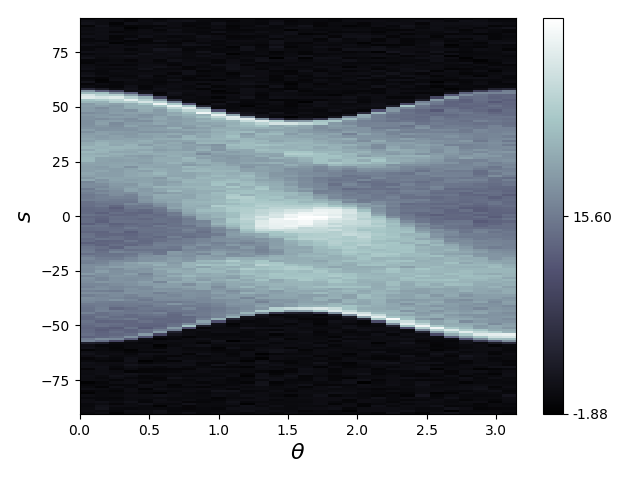}
			\caption{Shepp-Logan sinogram}
		\end{subfigure}
		\\
		\begin{subfigure}[t]{.48\linewidth}
			\includegraphics[width=\linewidth, trim={23mm 17mm 32mm 6mm}, clip]{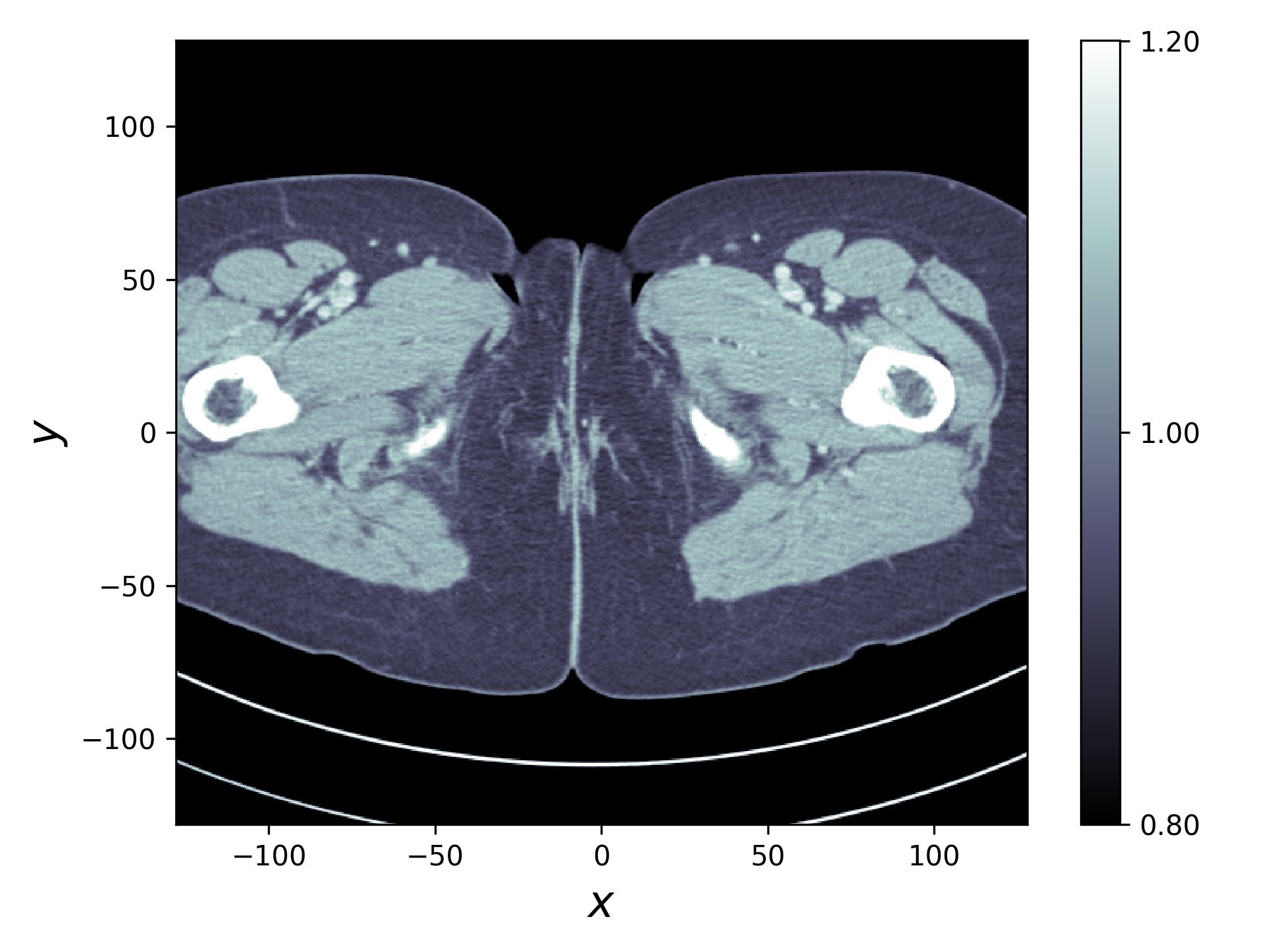}
			\caption{Human phantom}
		\end{subfigure}
		\begin{subfigure}[t]{.48\linewidth}	
			\includegraphics[width=\textwidth, trim={23mm 17mm 32mm 6mm}, clip]{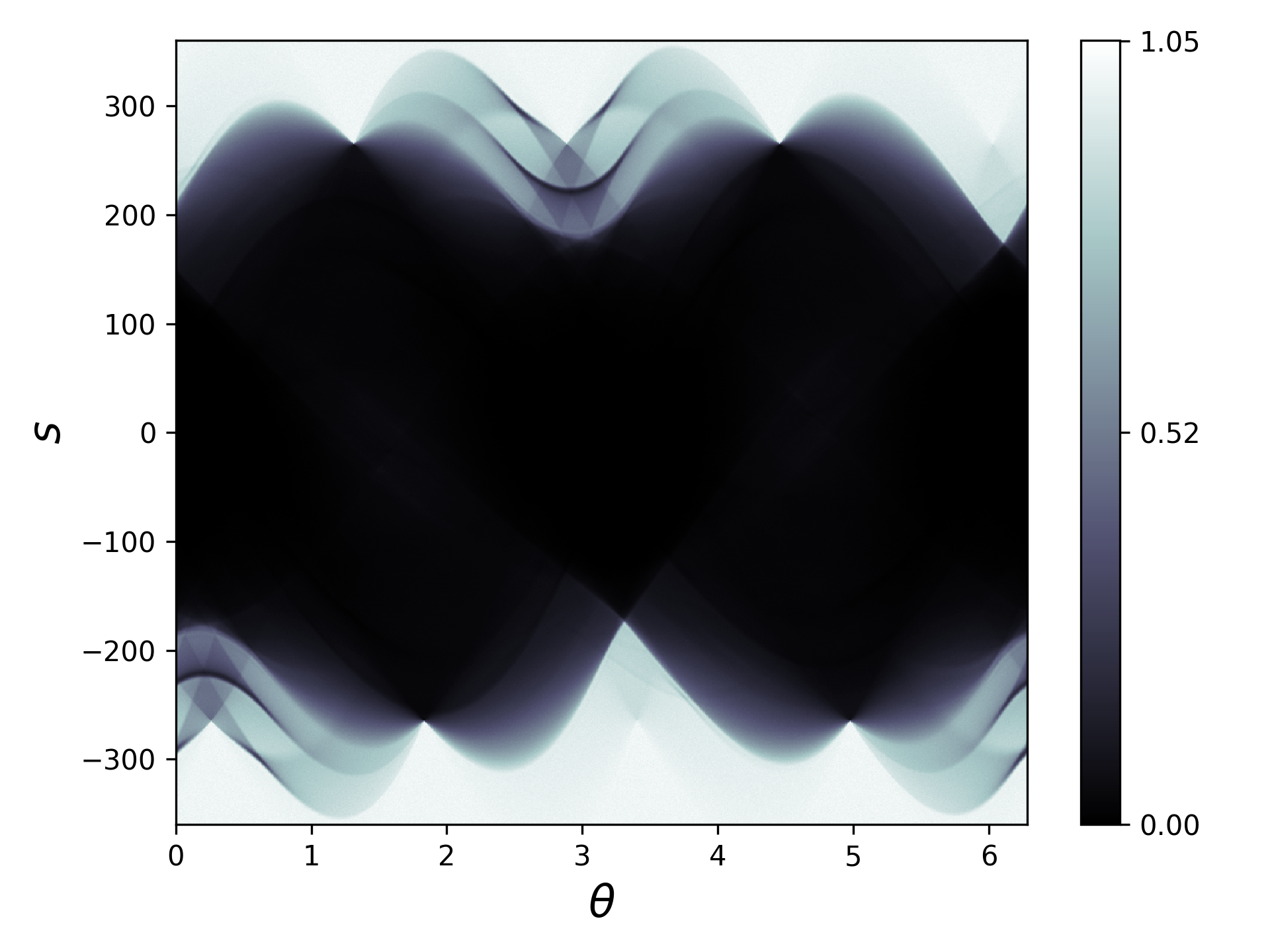}
			\caption{Human pre-log sinogram}
		\end{subfigure}
		\caption{Example of data which are used for training and validation. Top: Randomly generated ellipses. Middle: Validation data, here given by the modified Shepp-Logan phantom. Bottom: Validation data generated from the human phantoms.}
		\label{fig:example_data}
	\end{figure}

	We used a two-dimensional fan-beam geometry with 1000 angles, \unit{1000}{\pixels}, source to axis distance \unit{500}{\milli\meter} and axis to detector distance \unit{500}{\milli\meter}. In this setting, we consider the more physically correct non-linear forward model given by Beer-Lamberts law
	\begin{equation*}
	\label{eq:mayo_forward}
	\ForwardOp(\signal)(\ell) = e^{-\mu \RadonTransform(\signal)(\ell)}
	\end{equation*}
	where the unit of $f$ is \unit{}{\gram/\centi\cubic\meter} and $\mu$ is the mass attenuation coefficient, in this work selected to \unit{0.2}{\centi\squaren\meter/\gram} which is approximately the value in water at x-ray energies. We used Poisson noise corresponding to $10^4$ incident photons per pixel before attenuation, which would correspond to a low dose \ac{CT} scan. 
	We find the action of the adjoint of the derivative by straightforward computation
	\begin{equation*}
	\label{eq:mayo_backward}
	[\partial \ForwardOp(\signal)]^*(g) = - \mu \RadonTransform^* \bigl( e^{-\mu \RadonTransform(\signal)(\cdot)} \data(\cdot) \bigr) \quad \text{for $g \in \DataSpace$.}
	\end{equation*}
	The forward model can also be linearised by applying $-\log(\cdot) / \mu$ to both the data and forward operator, which then simply becomes the ray-transform as for the ellipse data. We implemented both the pre-log (non-linear) and post-log (linear) forward models and compare their results.	
	\\\\
	For validation of the ellipse data case, we simply use the (modified) Shepp-Logan phantom and for the human phantom data we use one held out set of patient data. See \cref{fig:example_data} for examples.
	
	\subsection{Implementation}\label{sec:Impl}
	The methods described above were implemented in Python using \ac{ODL} \cite{AdKoOk17} and TensorFlow \cite{AbAgBaBrCh15}. All operator-related components, such as the forward operator $\ForwardOp$, were implemented in \ac{ODL}, and these were then converted into TensorFlow layers using the \texttt{as\_tensorflow\_layer} functionality of \ac{ODL}. The neural network layers and training were implemented using TensorFlow. The implementation utilizes abstract \ac{ODL} structures for representing functional analytic notions and is therefore generic and easily adaptable to other inverse problems. In particular, the code can be easily adapted to other imaging modalities.
	
	We used the \ac{ODL} operator \texttt{RayTransform} in order to evaluate the ray transform and its adjoint using the \ac{GPU} accelerated \texttt{'astra\_gpu'} backend \cite{AaPaCaJaBl16}.

	\subsubsection{Deep neural network and training details}
	\label{sec:deep_neural_network}
	
	Given the general Learned Primal-Dual scheme in \cref{alg:learned_pd}, a parametrization of the learned proximal operators is needed in order to proceed. In many inverse problems, and particularly in \ac{CT} and \ac{MRI} reconstruction, most of the useful properties for both the forward operator and prior are approximately translation invariant. For this reason the resulting reconstruction operator should be approximately translation invariant, which indicates that \acp{CNN} are suitable for parametrizing the aforementioned reconstruction operator. 
	
	We used learned proximal operators of the form
	\[
	\IdentityOp + \AffineOp_{\weight_{3}, \bias_{3}} \circ \NonLinOp_{c_{2}} \circ \AffineOp_{\weight_{2}, \bias_{2}} \circ \NonLinOp_{c_{1}} \circ \AffineOp_{\weight_{1}, \bias_{1}}
	\]
	where $\IdentityOp$ is the identity operator that makes the network a \emph{residual} network. There are two main reasons for choosing such a structure. First, proximal operators (as the name implies) are typically close to the identity and second, there is rich evidence in the machine learning literature \cite{ResNet} that networks of this type are easier to train. Heuristically this is because each update does not need to learn the whole update, but only a small offset from the identity.
	
	Additionally, we used affine operators $\AffineOp_{\weight_{j}, \bias_{j}}$ parametrized by \emph{weights} $\weight_{j}$ and \emph{biases} $\bias_{j}$. The affine operators are defined in terms of so called convolution operators (here given on the primal space, but equivalently on the dual space). These are given as affine combinations of regular convolutions, more specifically:
	\[
	\AffineOp_{\weight_{j}, \bias_{j}} : \RecSpace^n \to \RecSpace^m 
	\]
	where the $k$:th component is given by
	\[ 
	\bigl[\AffineOp_{\weight_{j}, \bias_{j}}([\primal^{(1)}, \dots, \primal^{(n)}])\bigr]^{(k)}
	=
	\bias_{j}^{(k)} + \sum_{l = 1}^n \weight_{j}^{(l, k)} \ast \primal^{(l)}
	\]
	where $\bias_{j} \in \Real^m$ and $\weight_{j} \in \RecSpace^{n \times m}$.
	
	The non-linearities were chosen to be Parametric Rectified Linear Units (PReLU) functions
	\[
	\NonLinOp_{c_{j}}(x) = 
	\begin{cases}
	x & \text{if } x \geq 0 \\
	- c_{j} x & \text{else}.
	\end{cases}
	\]
	This type of non-linearity has proven successful in other applications such as classification \cite{PreluHE}.
	
	We let the number of data that persists between the iterates be $N_{\text{primal}} = N_{\text{dual}} = 5$. The convolutions were all $3 \times 3$ pixel size, and the number of channels was, for each primal learned proximal, $6 \to 32 \to 32 \to 5$, and for the duals $7 \to 32 \to 32 \to 5$ where the higher number of inputs is due to the data $g$ being supplied to the dual proximals.
	
\begin{figure*}
		\centering
		\includegraphics[width=\textwidth]{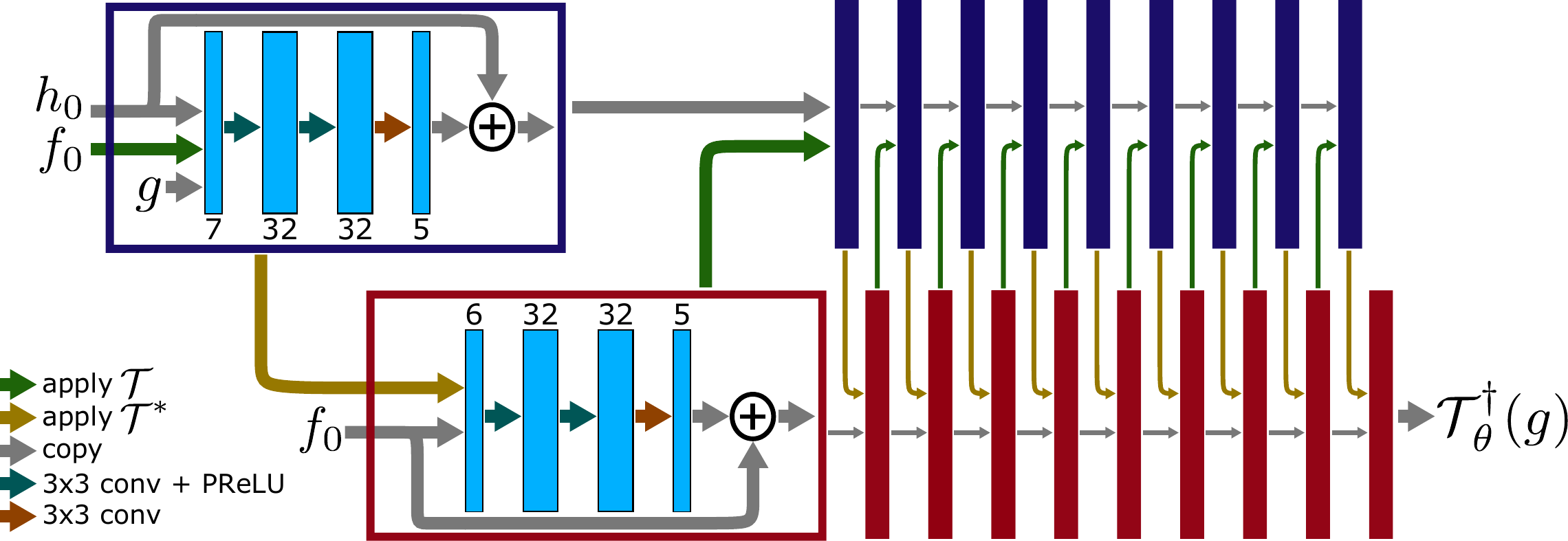}
		\caption{Network architecture used to solve the tomography problem. The dual iterates are given in blue boxes, while the primal iterates are in the red boxes. The blue/red boxes all have the same architecture, which is illustrated in the corresponding large boxes. Several arrows pointing to one box indicates concatenation. The initial guesses enter from the left, while the data is supplied to the dual iterates. In the classical \ac{PDHG} algorithm, the primal iterates would instead of a \ac{CNN} be given by $\ProxOp_{\tau \OpG}$  with over-relaxation, and the dual iterates would be given by $\ProxOp_{\sigma \OpF^*}$.}
		\label{fig:network_graph}
	\end{figure*}

	We let the number of unrolled iterations be $I=10$, that is the operator $\ForwardOp$ and the adjoint of its derivative $[\partial \ForwardOp(\primal_i^{(1)})]^*$ are both evaluated 10 times by the network. Since each iterate involves two 3-layer networks, one for each proximal, the total depth of the network is 60 convolutional layers and the total number of parameters approximately $2.4 \cdot 10^5$. In the context of deep learning, this is a deep network but with a small number of parameters. The network is visualized in \cref{fig:network_graph}.
	
	We used the Xavier initialization scheme \cite{XavierInit} for the convolution parameters, and initialized all biases to zero. 
	
	
	We trained the network by minimizing the empirical loss \cref{eq:LossEmpirical} using training data as explained above using the ADAM optimizer in TensorFlow \cite{ADAMOptimizer}. We used $10^5$ batches on each problem and used a learning rate schedule according to cosine annealing \cite{CosineAnnealing}, i.e. the learning rate in step $t$ was
	\[
	\eta_t = \frac{\eta_0}{2}\Bigl(1 + \cos\Bigl(\pi \frac{t}{t_{\max}} \Bigr) \Bigr)
	\]
	where the initial learning rate $\eta_0$ was set to $10^{-3}$. We also let the parameter $\beta_2$ of the ADAM optimizer to $0.99$ and let all other parameters use the default choices. We performed global gradient norm clipping \cite{GradientNormClip}, limiting the gradient norms to 1 in order to improve training stability and used a batch size of 5 for the ellipse data and 1 for the human phantoms.
	
	We did not use any regularization of the learned parameters, nor did we utilize any tricks such as dropout or batch normalization. Neither did we perform any data augmentation.
	
	The training was done using a single GTX 1080 TI \ac{GPU} and took about 11 hours for the ellipse data and 40 hours for the human phantoms.
	
	\subsubsection{Incorporating the forward operator in neural networks}	
	In order to minimize the loss function \cref{eq:LossStandard}, \ac{SGD} type methods are typically used and these require (an estimate of) the gradient of the loss function
	\[
	[\grad \loss](\vparam) =
	\Expect_{(\stsignal,\stdata) \sim \ProbabilityMeasure} 
	\Bigl[
	2 [\partial_{\vparam} \signallearned{\stdata}]^*\bigl( \signallearned{\stdata} - \stsignal \bigr)
	\Bigr],
	\]
	where $[\partial_{\vparam}\! \signallearned{\stdata}]^*$ is the adjoint of the derivative (with respect to $\vparam$) of the learned reconstruction operator applied in $g$, often called gradient in the machine learning literature. This introduces a challenge since it will depends on each component of the neural network, including the learned proximal operators but also the forward operator $\ForwardOp$ and the backward operator $[\partial \!\ForwardOp(\signal)]^*$, propagated through all $I$ iterations.
	
	To solve this, we used the built in automatic differentiation functionality of TensorFlow which uses the chain rule (back-propagation). This in turn requires the adjoints of the derivatives of each individual component which for the proximals were computed by TensorFlow and for the operators by ODL.

	\subsection{Comparison}
	
	We compare the algorithm to several widely used algorithms, including standard \ac{FBP} and (isotropic) \ac{TV} regularized reconstruction. We also compare against several learned schemes. These are briefly summarize here, see the references for full descriptions.
	
	The \ac{FBP} reconstructions were done with a Hann filter and used the method \texttt{fbp\_op} in ODL. The \ac{TV} reconstruction was performed using 1000 iterations of the classical \ac{PDHG} algorithm, implemented in ODL as \texttt{pdhg}. The filter bandwith in the \ac{FBP} reconstruction and the regularization parameter in the \ac{TV} reconstruction were selected in order to maximize the \ac{PSNR}. 
	
	The partially \emph{Learned Gradient} method in \cite{Adler17} is similar to the algorithm proposed in this article, but differs in that instead of learning proximal operators it learns a gradient operator and the forward operator enters into the neural network through the gradient of the data likelihood. Publicly available code and parameters \cite{SourceCodeGradient} were used.
	
	The next comparison is against a deep learning based approach for post-processing based on a so called \emph{U-Net} \cite{UNet}. The U-Net was first proposed for image segmentation, but by changing the number of output channels to one, it can also be used for post-processing as was done in \cite{JiMcFrUn16}.  Here an initial reconstruction is first performed using \ac{FBP} and a neural network is trained on pairs of noisy \ac{FBP} images and noiseless/low noise ground truth images, learning a mapping between them. We re-implemented the algorithm from \cite{JiMcFrUn16} but found that using the training procedure as stated in the paper gave sub-optimal results. We hence report values from using the same training scheme as for our other algorithms in order to give a more fair comparison.
	
	Additionally, our comparison includes learned \ac{PDHG}, \cref{alg:learned_cp}, as well as the following two simplified versions of the Learned Primal-Dual algorithm. The first is a \emph{Learned Primal} algorithm, which does not learn any parameters for the dual proximal, instead it returns the residual 
	\[
	\Gamma_{\vparam_i^d}\bigl(\dual_{i - 1}, \ForwardOp(\primal_{i-1}^{(2)}), \data\bigr) =\ForwardOp(\primal_{i-1}^{(2)}) - \data
	\]
	The second, \emph{FBP + residual denoising} algorithm, further simplifies the problem by discarding the forward operator completely, and can be seen as selecting
	\[
	\Lambda_{\vparam_i^p}\bigl(\primal_{i - 1}, [\partial\ForwardOp(\primal_{i-1}^{(1)})]^*(\dual_i^{(1)}) \bigr) =
	\Lambda_{\vparam_i^p}\bigl(\primal_{i - 1}\bigr)
	\]
	Since this method does not have access to the data $\data$, we select the initial guess according to a \ac{FBP}, see \cref{eq:init_pseudoinverse}. This makes the algorithm a learned denoiser.
	
	For the human phantoms we compare both non-linear and linearized versions of the forward operator, but given that training times are noticeably longer, we only compare to the previously established methods of \ac{FBP}, \ac{TV} and U-Net denoising.
	
	All learned algorithms were trained using the same training scheme as outlined in \cref{sec:deep_neural_network}, and measure the run-time, \ac{PSNR} and the \ac{SSIM} \cite{SSIM}.
	
	All methods that we compare against are available in the accompanying source code.
	
	\section{Results}
	
	The quantitative results for the ellipse data is given in \cref{tab:quant_res_rt}, where we can see that the proposed Learned Primal-Dual algorithm out-performs the classical schemes (\ac{FBP} and \ac{TV}) significantly w.r.t. the reconstruction error as measured by both \ac{PSNR} and \ac{SSIM}. We also note that the Learned Primal-Dual algorithm gives a large improvement over the previous deep learning based methods such as the learned gradient scheme and U-Net based post-processing, giving an improvement exceeding \unit{6}{\decibel}. The Learned Primal-Dual algorithm also outperforms the Learned \ac{PDHG} and the FBP + residual denoising algorithms by wide margins.
	
	The only method that achieves results close to the Learned Primal-Dual algorithm is the Learned Primal method, but the Learned Primal-Dual algorithm gives a noticeable improvement of \unit{1.3}{\decibel}.
	
\begin{figure}
		\centering	
		\begin{subfigure}[t]{.48\linewidth}
			\includegraphics[width=\linewidth, trim={23mm 17mm 32mm 6mm}, clip]{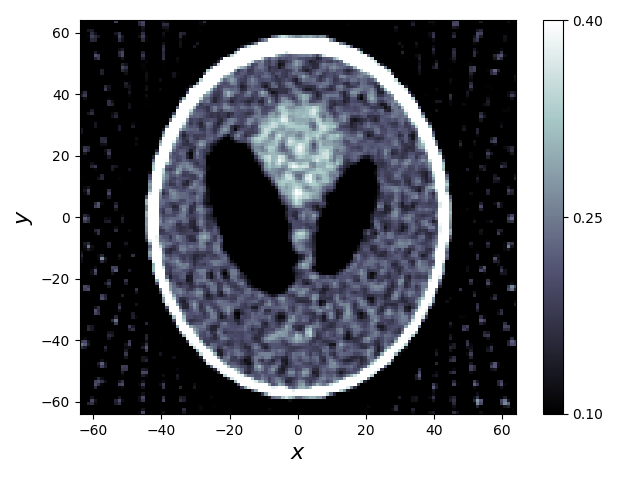}
			\caption{\Ac{FBP}}
			\label{fig:shepp_logan_fbp_windowed}
		\end{subfigure}
		\begin{subfigure}[t]{.48\linewidth}	
			\includegraphics[width=\textwidth, trim={23mm 17mm 32mm 6mm}, clip]{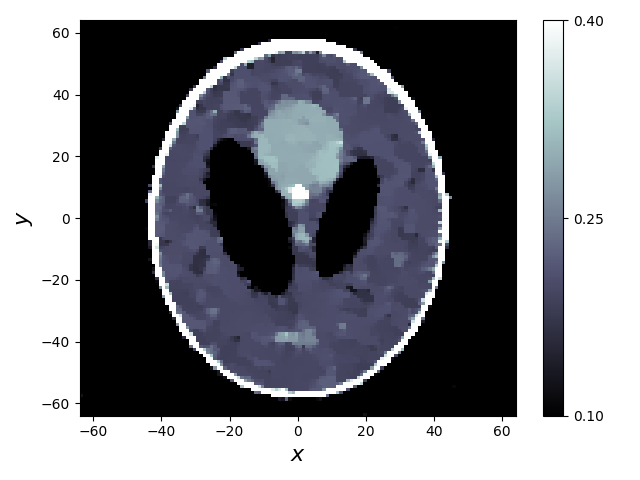}
			\caption{\Ac{TV}}
			\label{fig:shepp_logan_tv_windowed}
		\end{subfigure}
		\\
		\begin{subfigure}[t]{.48\linewidth}	
			\includegraphics[width=\linewidth, trim={23mm 17mm 32mm 6mm}, clip]{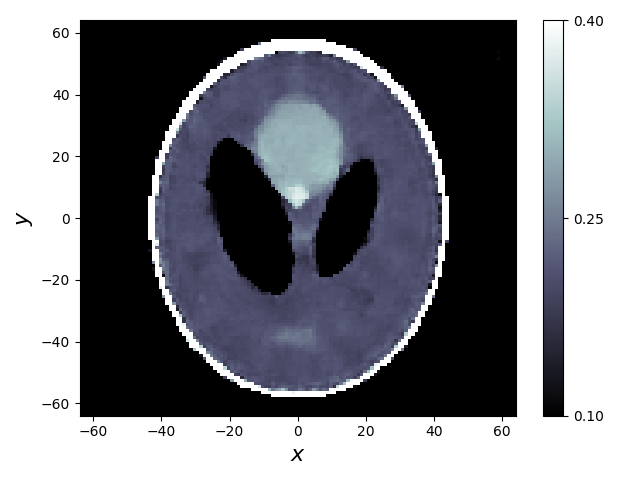}
			\caption{FBP + U-Net denoising}
			\label{fig:shepp_logan_result_u_net}
		\end{subfigure}
		\begin{subfigure}[t]{.48\linewidth}	
			\includegraphics[width=\linewidth, trim={23mm 17mm 32mm 6mm}, clip]{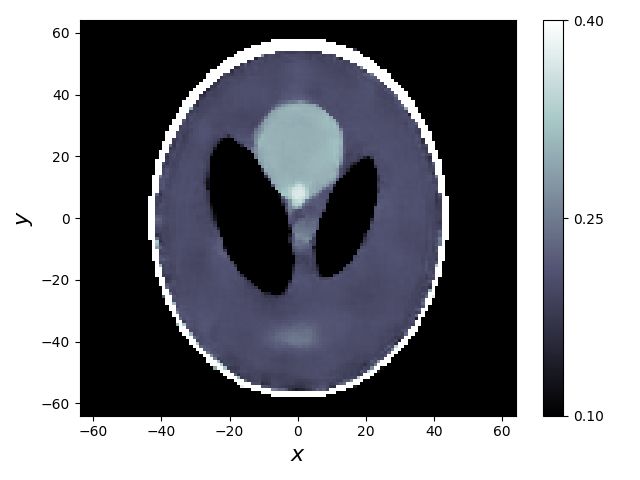}
			\caption{FBP + residual denoising}
			\label{fig:primal_noop_ellipses_result}
		\end{subfigure}
		\\
		\begin{subfigure}[t]{.48\linewidth}
			\includegraphics[width=\textwidth, trim={23mm 17mm 32mm 6mm}, clip]{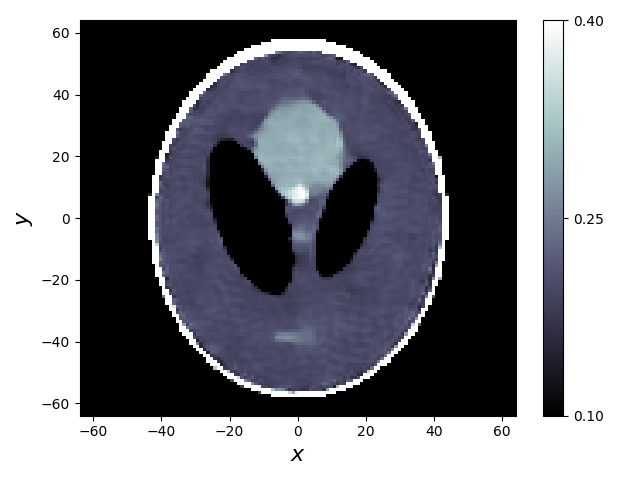}
			\caption{Learned gradient}
			\label{fig:shepp_logan_result_windowed_iterate_10}
		\end{subfigure}
		\begin{subfigure}[t]{.48\linewidth}	
			\includegraphics[width=\linewidth, trim={23mm 17mm 32mm 6mm}, clip]{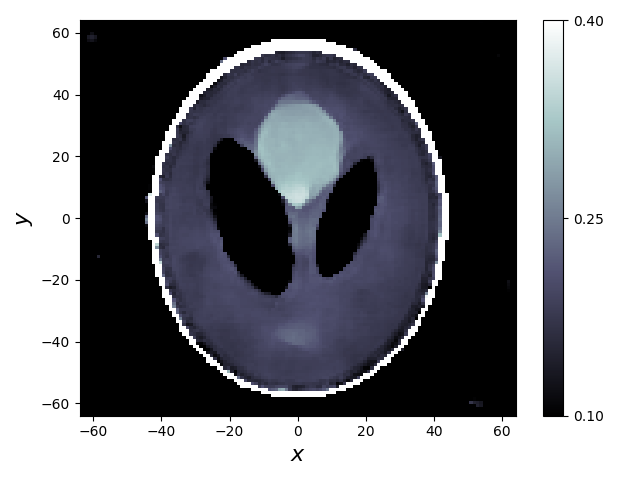}
			\caption{Learned \ac{PDHG}}
			\label{fig:shepp_logan_learned_cp}
		\end{subfigure}
		\\
		\begin{subfigure}[t]{.48\linewidth}
			\includegraphics[width=\textwidth, trim={23mm 17mm 32mm 6mm}, clip]{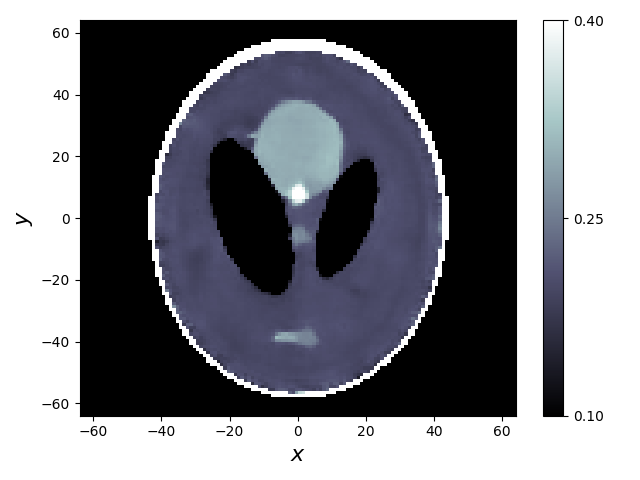}
			\caption{Learned Primal}
			\label{fig:shepp_logan_learned_primal}
		\end{subfigure}
		\begin{subfigure}[t]{.48\linewidth}
			\includegraphics[width=\textwidth, trim={23mm 17mm 32mm 6mm}, clip]{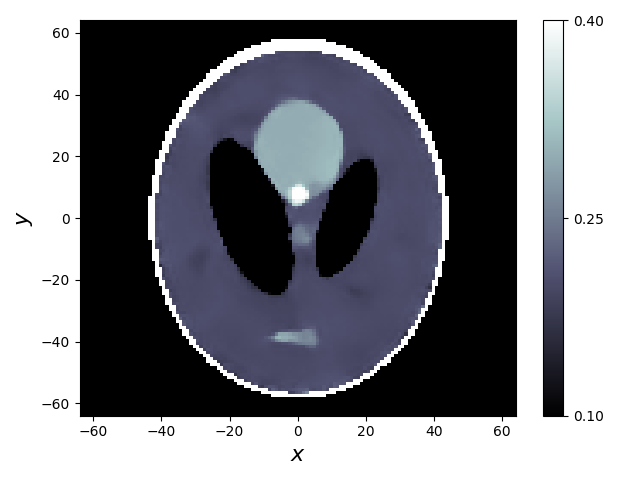}
			\caption{Learned Primal-Dual}
			\label{fig:primal_dual_ellipses_result}
		\end{subfigure}
		\caption{Reconstructions for the ellipse data using the compared methods. The window is set to $[0.1, 0.4]$, corresponding to the soft tissue of the modified Shepp-Logan phantom.}
		\label{fig:ray_transform_results}
	\end{figure}

	The results are visualized in \cref{fig:ray_transform_results}. We note that small structures, such as the small inserts, are much more clearly visible in the Learned Primal and Learned Primal-Dual reconstructions than in the other reconstructions. We also note that both the Learned \ac{PDHG} and Learned Primal reconstruction seem to have a halo artefact close to the outer bone which is absent in the Learned Primal-Dual reconstruction.
	
	With respect to run-time the learned methods that involve calls to the forward operator (Learned Gradient, \ac{PDHG}, Primal, Primal-Dual) are slower than the methods that do not (FBP + U-Net denoising, Residual) by a factor $\approx 6$. When compared to \ac{TV} regularized reconstruction all learned methods are at least 2 orders of magnitude faster.
	
	\begin{table}
		\centering
		\npdecimalsign{.}
		\nprounddigits{3}
		\begin{tabular}{ l r n{1}{3} r r}
			\toprule
			Method & \acs{PSNR} & SSIM & Runtime & Parameters\\
			\midrule
			\ac{FBP}  & 19.75 & 0.597091229970 & \textbf{4} & \textbf{1} \\
			\ac{TV} & 28.06 & 0.9289187090364599 & 5\,166 & \textbf{1} \\
			FBP + U-Net denoising & 29.20 & 0.943696335723 & 9 & $10^7$\\
			FBP + residual denoising & 32.38 & 0.971950173596 & 9 & $1.2 \cdot 10^5$ \\
			Learned Gradient & 32.29 & 0.981379635523 & 56 & $1.2 \cdot 10^4$\\
			Learned \ac{PDHG} & 28.32 & 0.909222619552 & 48 & $2.4 \cdot 10^4$\\
			Learned Primal & 36.97 & 0.986370310292 & 43 & $1.2 \cdot 10^5$\\
			Learned Primal-Dual\phantom{, linear} & \textbf{38.28} & {\npboldmath}0.988749956718 & 49 & $2.4 \cdot 10^5$\\
			\bottomrule
		\end{tabular}
		\caption{Comparison of reconstruction methods for the ellipses. \Ac{PSNR} measured in {\decibel} and runtime in \milli\second.}
		\label{tab:quant_res_rt}
	\end{table}
	
	\begin{table}
		\centering
		\npdecimalsign{.}
		\nprounddigits{3}
		\setlength\tabcolsep{5.0pt} 
		\begin{tabular}{ l r n{1}{3} r r}
			\toprule
			Method & \acs{PSNR} & SSIM & Runtime & Parameters\\
			\midrule
			\ac{FBP}  & 33.65 & 0.829697955203 & \textbf{423} & \textbf{1} \\
			\ac{TV} & 37.48 & 0.9463029565752 & 64\,371 & \textbf{1}\\
			FBP + U-Net denoising & 41.92 & 0.9413786 & 463 & $10^7$\\
			Learned Primal-Dual, linear & \textbf{44.11} & {\npboldmath}0.9694 & 620 & $2.4 \cdot 10^5$\\
			Learned Primal-Dual, non-linear & 43.91 & {\npboldmath}0.96873 & 670 & $2.4 \cdot 10^5$\\
			\bottomrule
		\end{tabular}
		\caption{Comparison of the Learned Primal-Dual algorithm with other methods for the Human phantom data. Units for entries are the same as in \cref{tab:quant_res_rt}.}
		\label{tab:quant_res_mayo}
	\end{table}
	
	\begin{figure*}
		\centering	
		\begin{subfigure}[t]{0.328\linewidth}
			\centering	
			\begin{tikzpicture}[
			remember picture,
			spy using outlines={%
				circle,
				red,
				magnification=4,
				size=1.5cm,
				connect spies,
				spy connection path={\draw[thick] (tikzspyonnode) -- (tikzspyinnode);}
			}
			]
			\node {\includegraphics[width=\linewidth, trim={22.5mm 17mm 27mm 6mm}, clip]{figures/mayo_phantom}};
			\spy on (-1.44,0.62) in node [left] at (-1.0,2.2);
			\spy on (1.69,1.04) in node [left] at (2.4,2.2);
			\end{tikzpicture}
			\caption{$\unit{512 \times 512}{\pixel}$ human phantom}
		\end{subfigure}
		\begin{subfigure}[t]{0.328\linewidth}
			\centering		
			\begin{tikzpicture}[
			remember picture,
			spy using outlines={%
				circle,
				red,
				magnification=4,
				size=1.5cm,
				connect spies,
				spy connection path={\draw[thick] (tikzspyonnode) -- (tikzspyinnode);}
			}
			]
			\node {\includegraphics[width=\linewidth, trim={22.5mm 17mm 27mm 6mm}, clip]{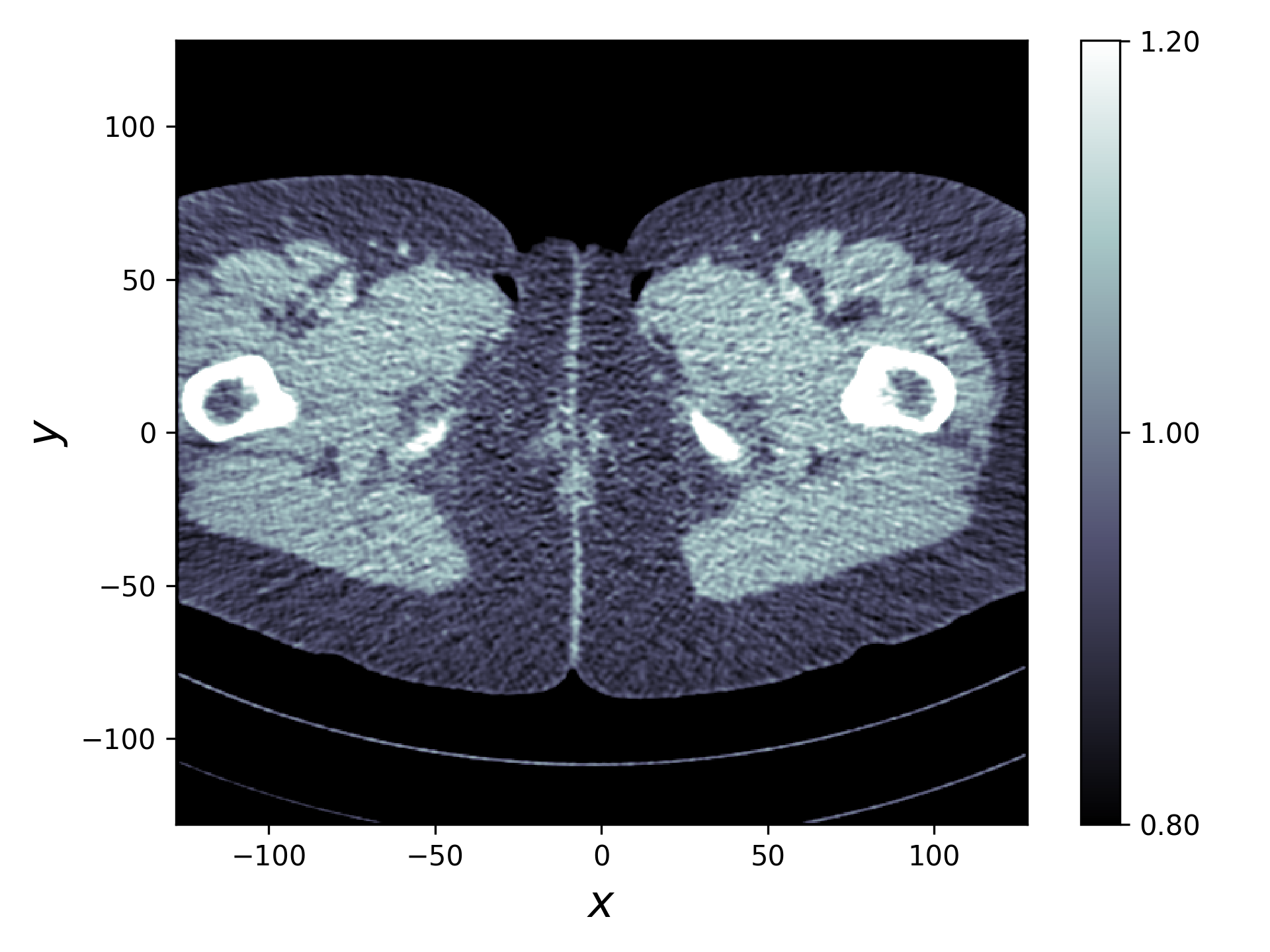}};
			\spy on (-1.44,0.62) in node [left] at (-1.0,2.2);
			\spy on (1.69,1.04) in node [left] at (2.4,2.2);
			\end{tikzpicture}
			\caption{\Acf{FBP}\\ \ac{PSNR} \unit{33.65}{\decibel}, \ac{SSIM} 0.830, \unit{423}{\milli\second}}
		\end{subfigure}
		\begin{subfigure}[t]{0.328\linewidth}
			\centering
			\begin{tikzpicture}[
			remember picture,
			spy using outlines={%
				circle,
				red,
				magnification=4,
				size=1.5cm,
				connect spies,
				spy connection path={\draw[thick] (tikzspyonnode) -- (tikzspyinnode);}
			}
			]
			\node {\includegraphics[width=\linewidth, trim={22.5mm 17mm 27mm 6mm}, clip]{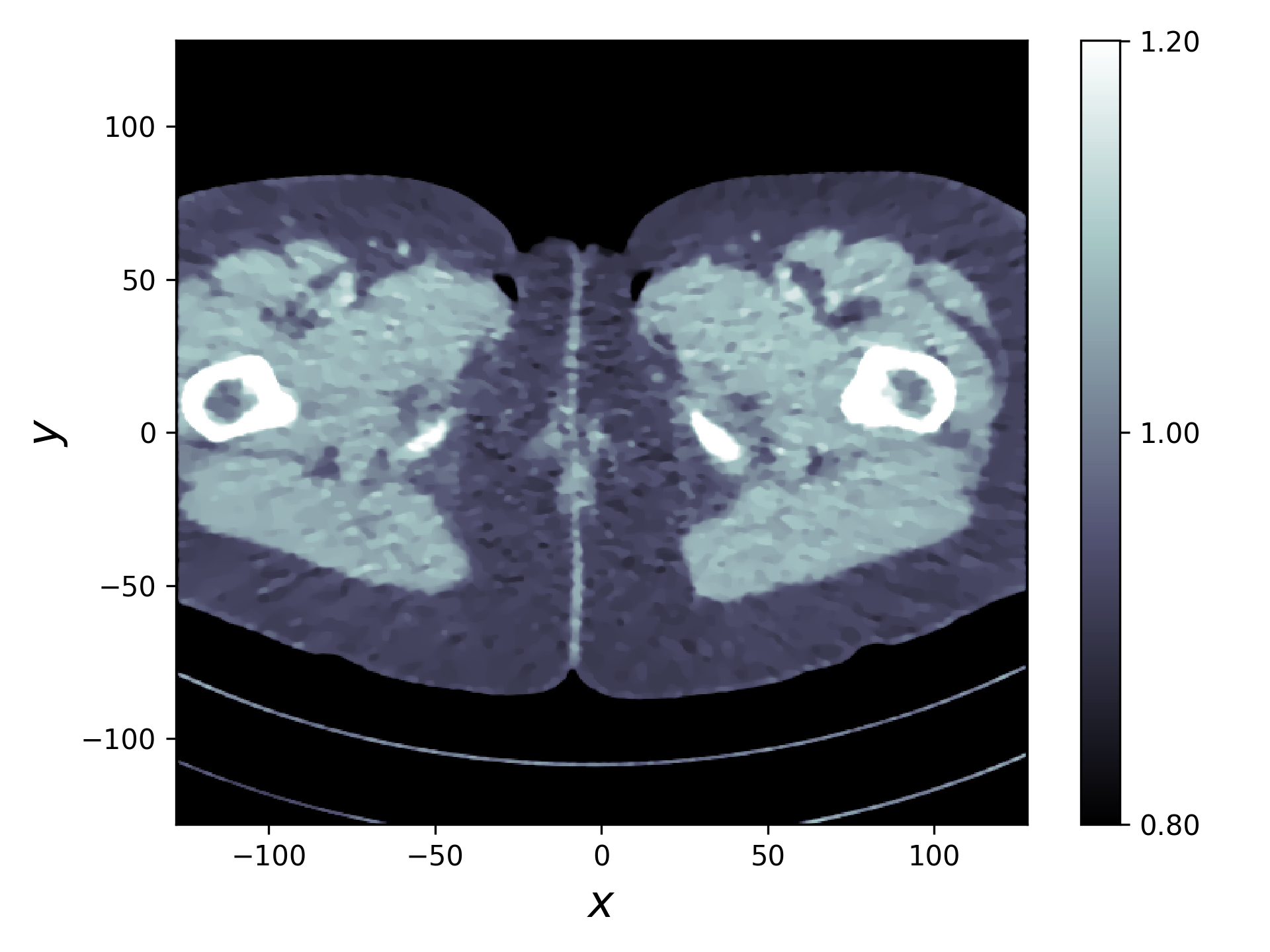}};
			\spy on (-1.44,0.62) in node [left] at (-1.0,2.2);
			\spy on (1.69,1.04) in node [left] at (2.4,2.2);
			\end{tikzpicture}
			\caption{\Acf{TV}\\ \ac{PSNR} \unit{37.48}{\decibel}, \ac{SSIM} 0.946, \unit{64\,371}{\milli\second}}
		\end{subfigure}
		\\[0.75em]
		\begin{subfigure}[t]{0.328\linewidth}
			\centering	
			\begin{tikzpicture}[
			remember picture,
			spy using outlines={%
				circle,
				red,
				magnification=4,
				size=1.5cm,
				connect spies,
				spy connection path={\draw[thick] (tikzspyonnode) -- (tikzspyinnode);}
			}
			]
			\node {\includegraphics[width=\linewidth, trim={22.5mm 17mm 27mm 6mm}, clip]{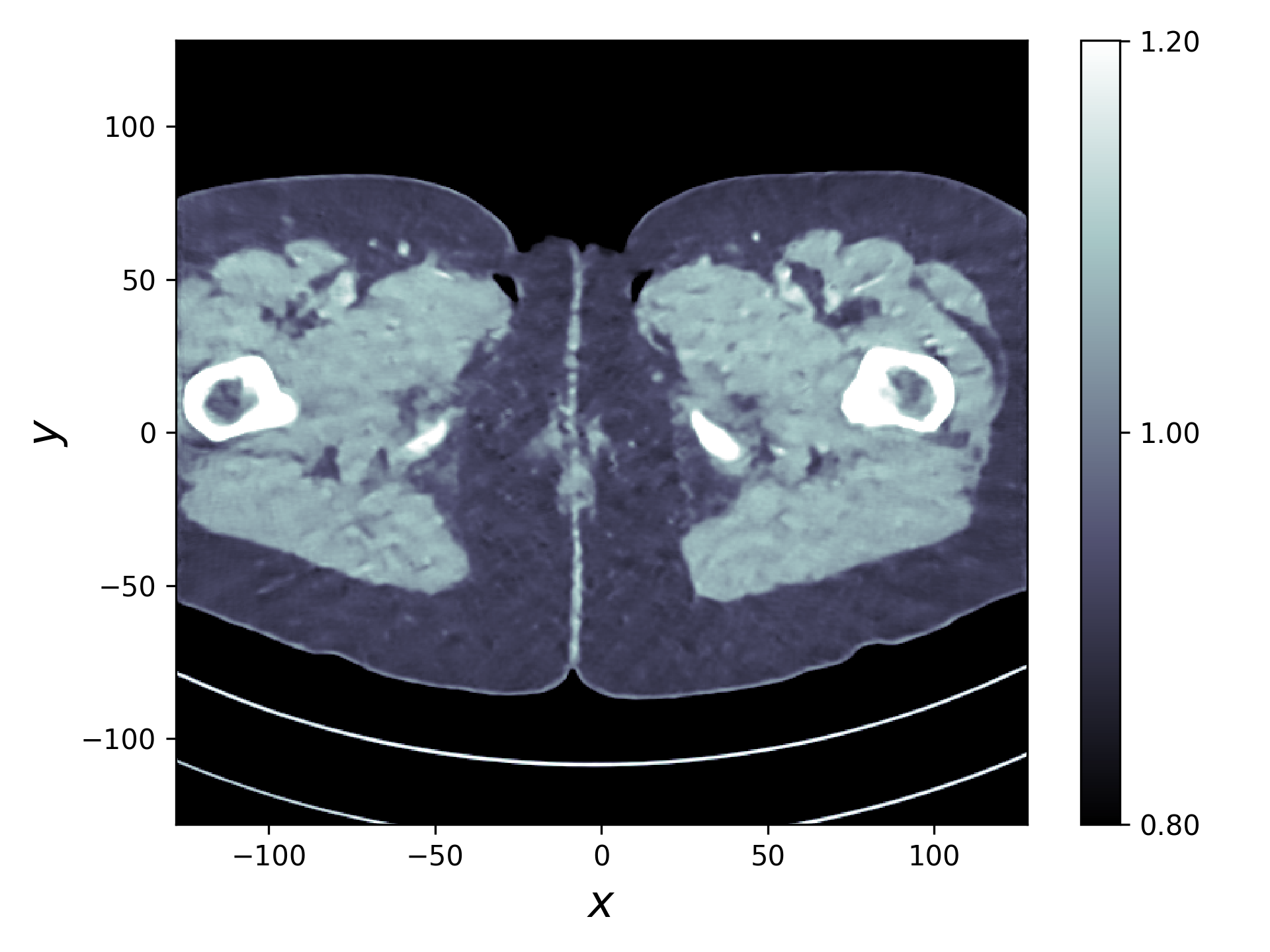}};
			\spy on (-1.44,0.62) in node [left] at (-1.0,2.2);
			\spy on (1.69,1.04) in node [left] at (2.4,2.2);
			\end{tikzpicture}
			\caption{FBP + U-Net denoising\\ \ac{PSNR} \unit{41.92}{\decibel}, \ac{SSIM} 0.941, \unit{463}{\milli\second}}
		\end{subfigure}
		\begin{subfigure}[t]{0.328\linewidth}
			\centering
			\begin{tikzpicture}[
			remember picture,
			spy using outlines={%
				circle,
				red,
				magnification=4,
				size=1.5cm,
				connect spies,
				spy connection path={\draw[thick] (tikzspyonnode) -- (tikzspyinnode);}
			}
			]
			\node {\includegraphics[width=\linewidth, trim={22.5mm 17mm 27mm 6mm}, clip]{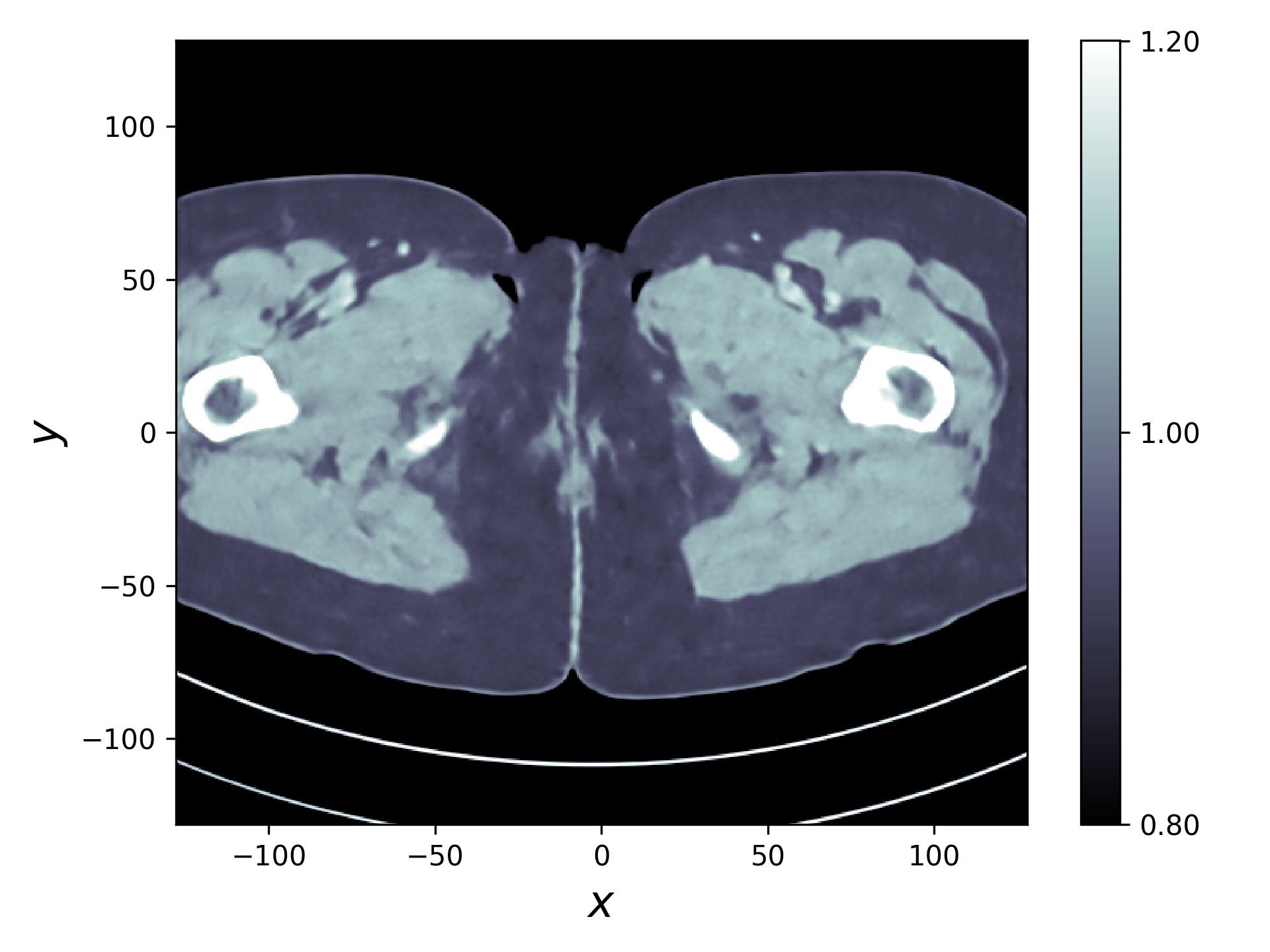}};
			\spy on (-1.44,0.62) in node [left] at (-1.0,2.2);
			\spy on (1.69,1.04) in node [left] at (2.4,2.2);
			\end{tikzpicture}
			\caption{Primal-Dual, linear\\ \ac{PSNR} \unit{44.10}{\decibel}, \ac{SSIM} 0.969, \unit{620}{\milli\second}}
		\end{subfigure}
		\begin{subfigure}[t]{0.328\linewidth}
			\centering
			\begin{tikzpicture}[
			remember picture,
			spy using outlines={%
				circle,
				red,
				magnification=4,
				size=1.5cm,
				connect spies,
				spy connection path={\draw[thick] (tikzspyonnode) -- (tikzspyinnode);}
			}
			]
			\node {\includegraphics[width=\linewidth, trim={22.5mm 17mm 27mm 6mm}, clip]{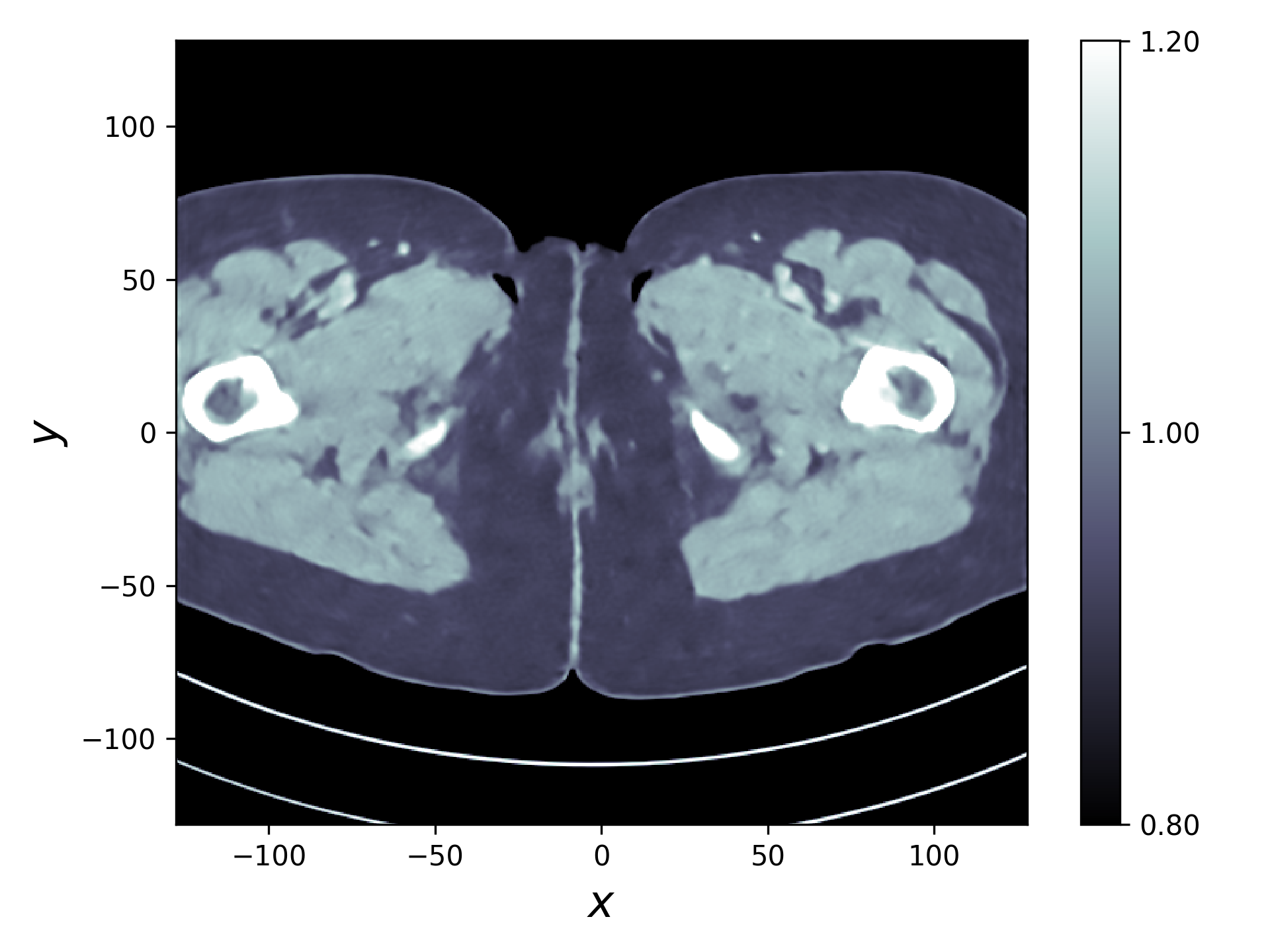}};
			\spy on (-1.44,0.62) in node [left] at (-1.0,2.2);
			\spy on (1.69,1.04) in node [left] at (2.4,2.2);
			\end{tikzpicture}
			\caption{Primal-Dual, non-linear\\ \ac{PSNR} \unit{43.91}{\decibel}, \ac{SSIM} 0.969, \unit{670}{\milli\second}}
		\end{subfigure}
		\caption{Reconstructions of a human phantom along with two zoomed in regions indicated by small circles. The left zoom-in has a true feature whereas texture in right zoom-in is uniform. The window is set to $\unit{[-200, 200]}{\hounsfield}$. Among the methods tested, only the Learned Primal-Dual algorithm correctly recovers these regions. In the others, the true feature in the left zoom-in is indistinguishable from other false features of same size/contrast and right-zoom in has a streak artifact. Note also the clinically feasible runtime of the Learned Primal-Dual algorithm. To summarize, the Learned Primal-Dual algorithm offers performance advantages over other methods that translate into true clinical usefulness.
		}
		\label{fig:results_mayo}
	\end{figure*}
	
	Quantitative results for the human phantoms data are presented in \cref{tab:quant_res_mayo}. We note that the \ac{FBP} reconstruction has a much more competitive image quality than it had for the ellipse data, both quantitatively and visually. It is likely for this reason that the FBP + U-Net denoising performs better than it did on the ellipses, outperforming \ac{TV} by \unit{4.4}{\decibel}. However, if we look at the \ac{SSIM} we note that this improvement does not translate as well to the structural similarity, where the method is comparable to \ac{TV} regularization.
	
	Both quantitatively and visually, the linear and non-linear versions of the Learned Primal-Dual algorithm give very similar results. We will focus on the linear version which gave slightly better results.
	
	The Learned Primal-Dual algorithm gives a \unit{10.5}{\decibel} improvement over the \ac{FBP} reconstruction, a \unit{6.6}{\decibel} improvement over \ac{TV} and \unit{2.2}{\decibel} over the U-Net denoiser. This is less than for the ellipse data, but still represents a large improvement. On the other hand, while the U-Net denoiser did not improve the \ac{SSIM} as compared to \ac{TV} regularization, the Learned Primal-Dual algorithm gives a large improvement.
	
	This improvement is also present in the images when inspected visually in \cref{fig:results_mayo}. In particular, we see that some artifacts visible in the \ac{FBP} reconstruction are still discernible in the U-Net denoiser and \ac{TV} reconstructions. Examples include streak artifacts, especially around the edges of the phantom and structures spuriously created from noise, such as a line in the muscle above the right bone. These are mostly absent in the Learned Primal-Dual reconstruction. However, we do note that the images do look slightly over-smoothed. Both of these observations become especially apparent if we look at the zoomed in regions, where we note that the Learned Primal-Dual algorithm is able to reconstruct finer detail than the other algorithms, but gives a very smooth texture.
	
	With respect to the run time, the Learned Primal-Dual is more competitive with the \ac{FBP} and U-Net denoiser algorithms for full size data than for the ellipse data. This is because the size of the data is much larger, which increases the runtime of the \ac{FBP} reconstruction, which is also needed to compute the initial guess for the U-Net denoiser. As for the ellipse data, both learned methods outperform \ac{TV} regularized reconstruction by two orders of magnitude with respect to runtime.
	
	\section{Discussion}
	The results show that the Learned Primal-Dual algorithm outperforms classical reconstruction algorithm by large margins as measured in both \ac{PSNR} and \ac{SSIM} and also improves upon learned post-processing methods for both simplified ellipse data and for human phantoms. In addition, especially for the $512 \times 512$ human phantoms, the reconstruction time is comparable with even filtered back-projection and learned post-processing.
	
	One interesting, and to the best of our knowledge, unique feature of the Learned Primal-Dual algorithm in the field of deep learning based \ac{CT} reconstruction, is that it gives reconstructions working directly from data, without any initial reconstruction as input.
	
\begin{figure}
		\centering
		\begin{subfigure}[t]{.32\linewidth}	
			\includegraphics[width=\linewidth, trim={23mm 17mm 32mm 6mm}, clip]{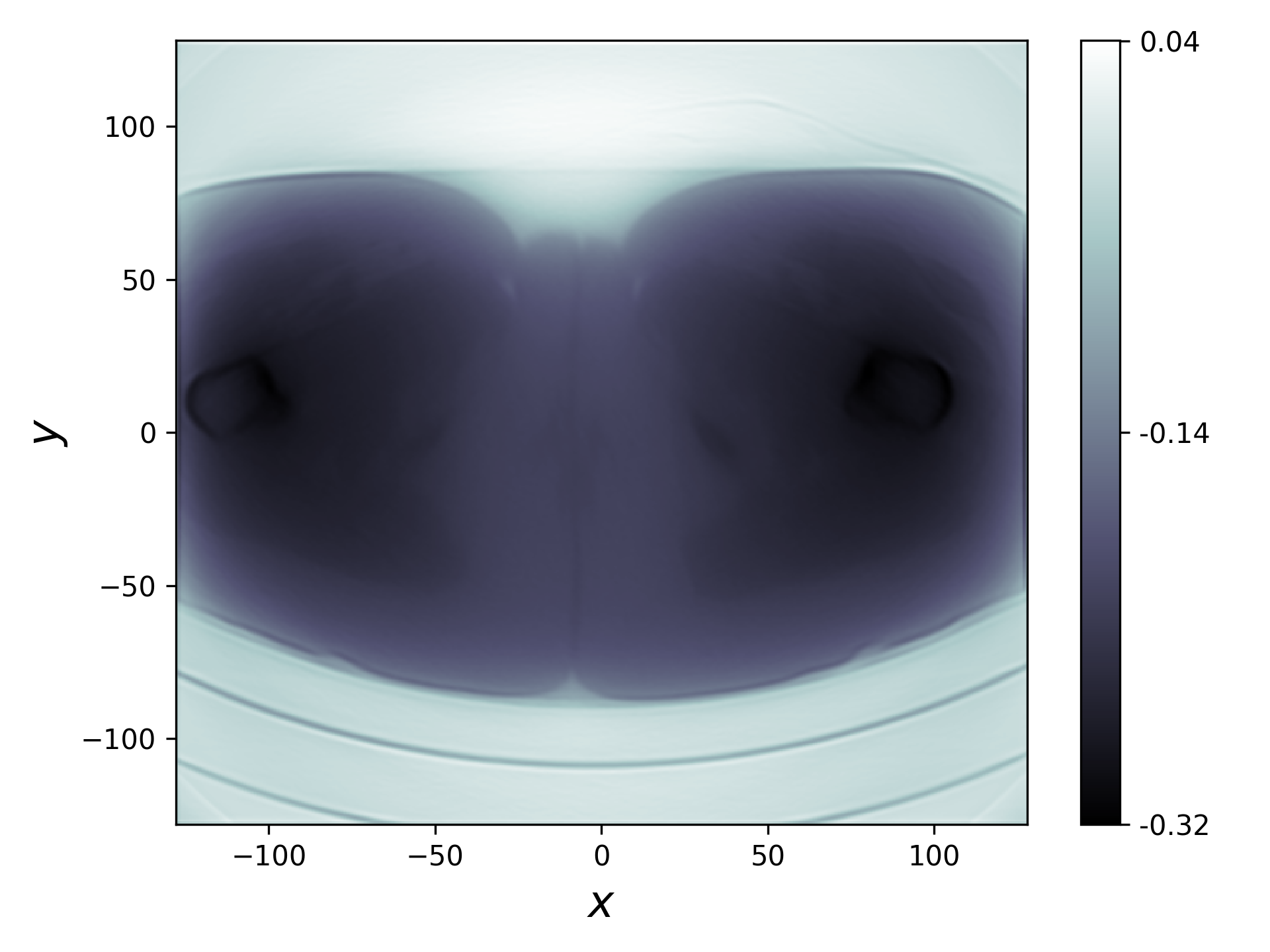}
		\end{subfigure}
		\begin{subfigure}[t]{.32\linewidth}
			\includegraphics[width=\textwidth, trim={23mm 17mm 32mm 6mm}, clip]{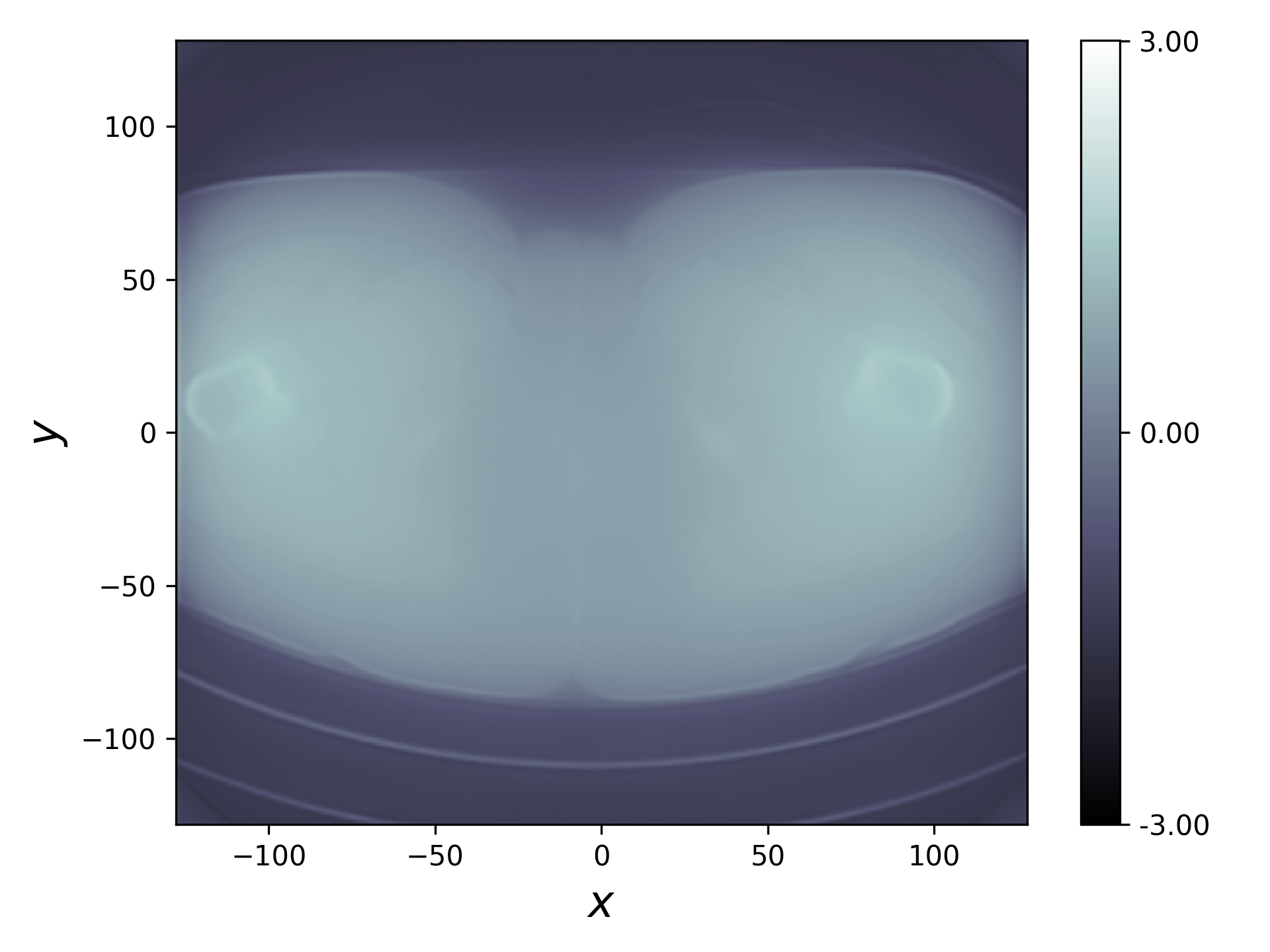}
		\end{subfigure}
		\begin{subfigure}[t]{.32\linewidth}
			\includegraphics[width=\textwidth, trim={23mm 17mm 32mm 6mm}, clip]{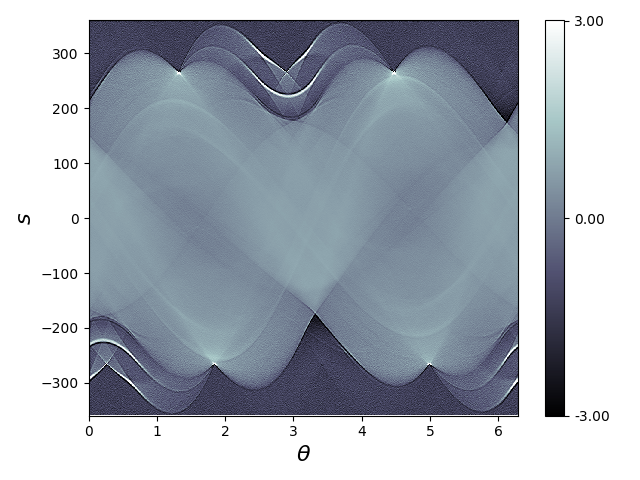}
		\end{subfigure}
		\\ \vspace{1mm}
		\begin{subfigure}[t]{.32\linewidth}	
			\includegraphics[width=\linewidth, trim={23mm 17mm 32mm 6mm}, clip]{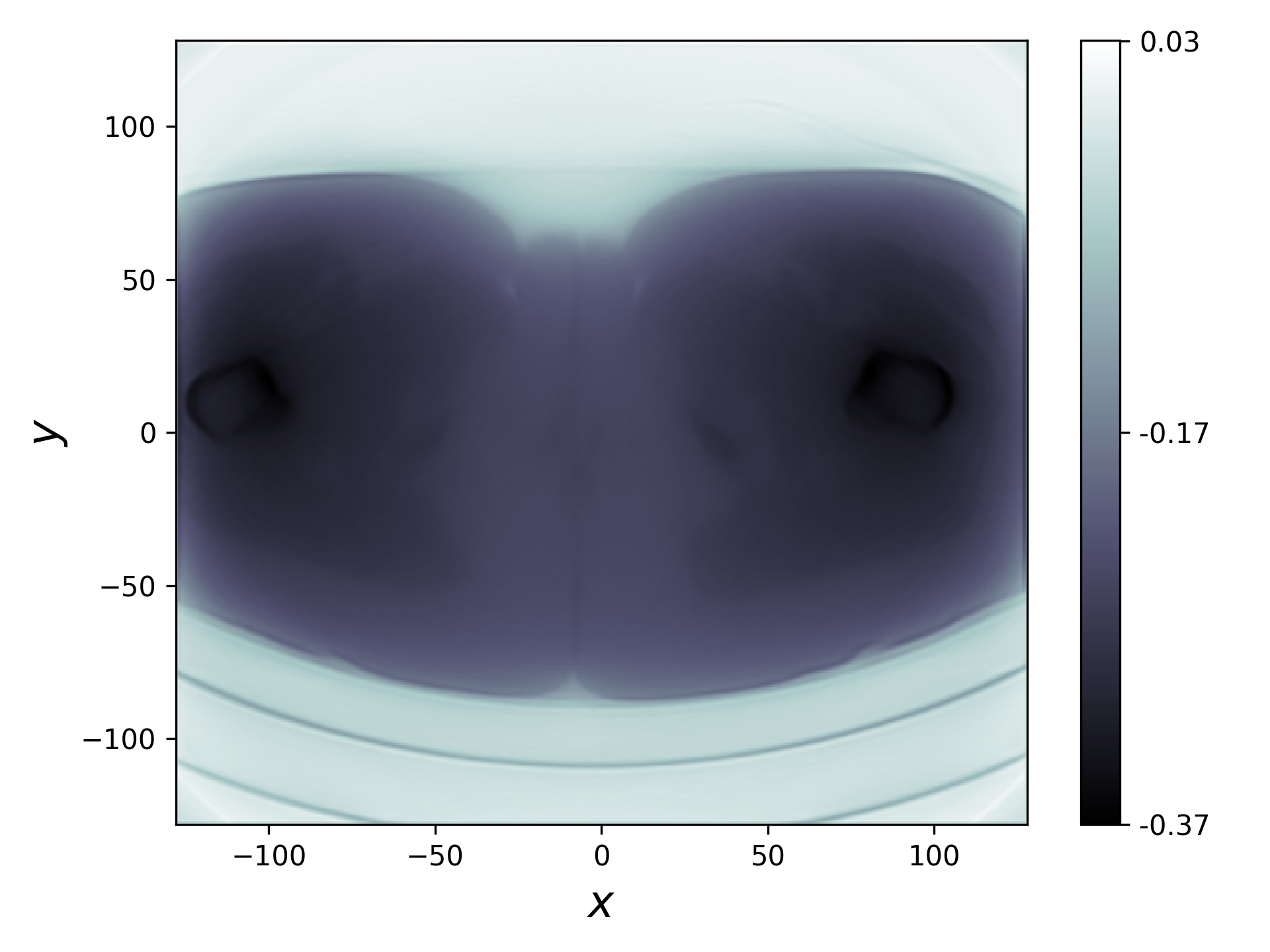}
		\end{subfigure}
		\begin{subfigure}[t]{.32\linewidth}
			\includegraphics[width=\textwidth, trim={23mm 17mm 32mm 6mm}, clip]{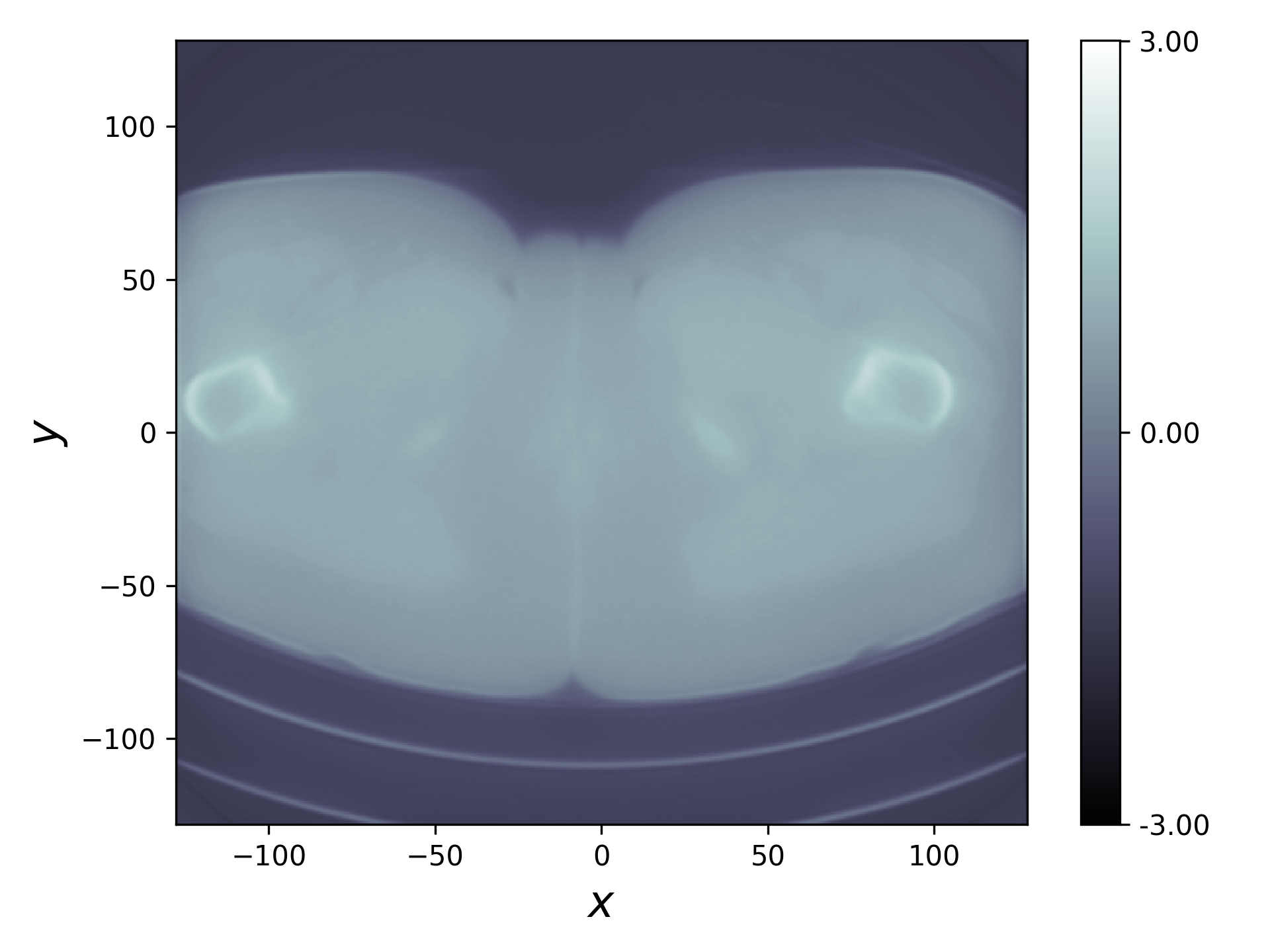}
		\end{subfigure}
		\begin{subfigure}[t]{.32\linewidth}
			\includegraphics[width=\textwidth, trim={23mm 17mm 32mm 6mm}, clip]{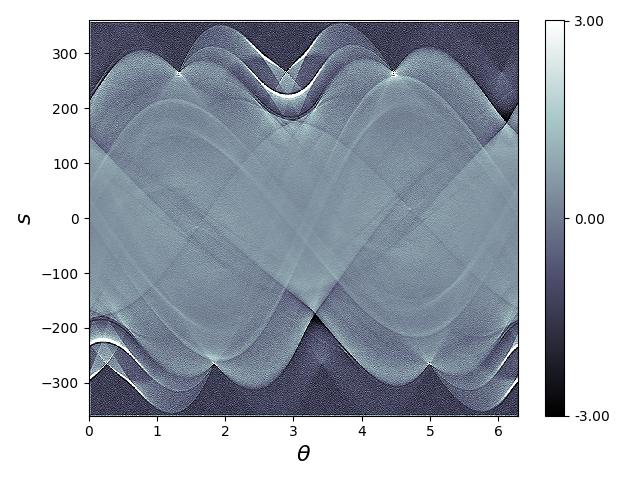}
		\end{subfigure}
		\\ \vspace{1mm}
		\begin{subfigure}[t]{.32\linewidth}	
			\includegraphics[width=\linewidth, trim={23mm 17mm 32mm 6mm}, clip]{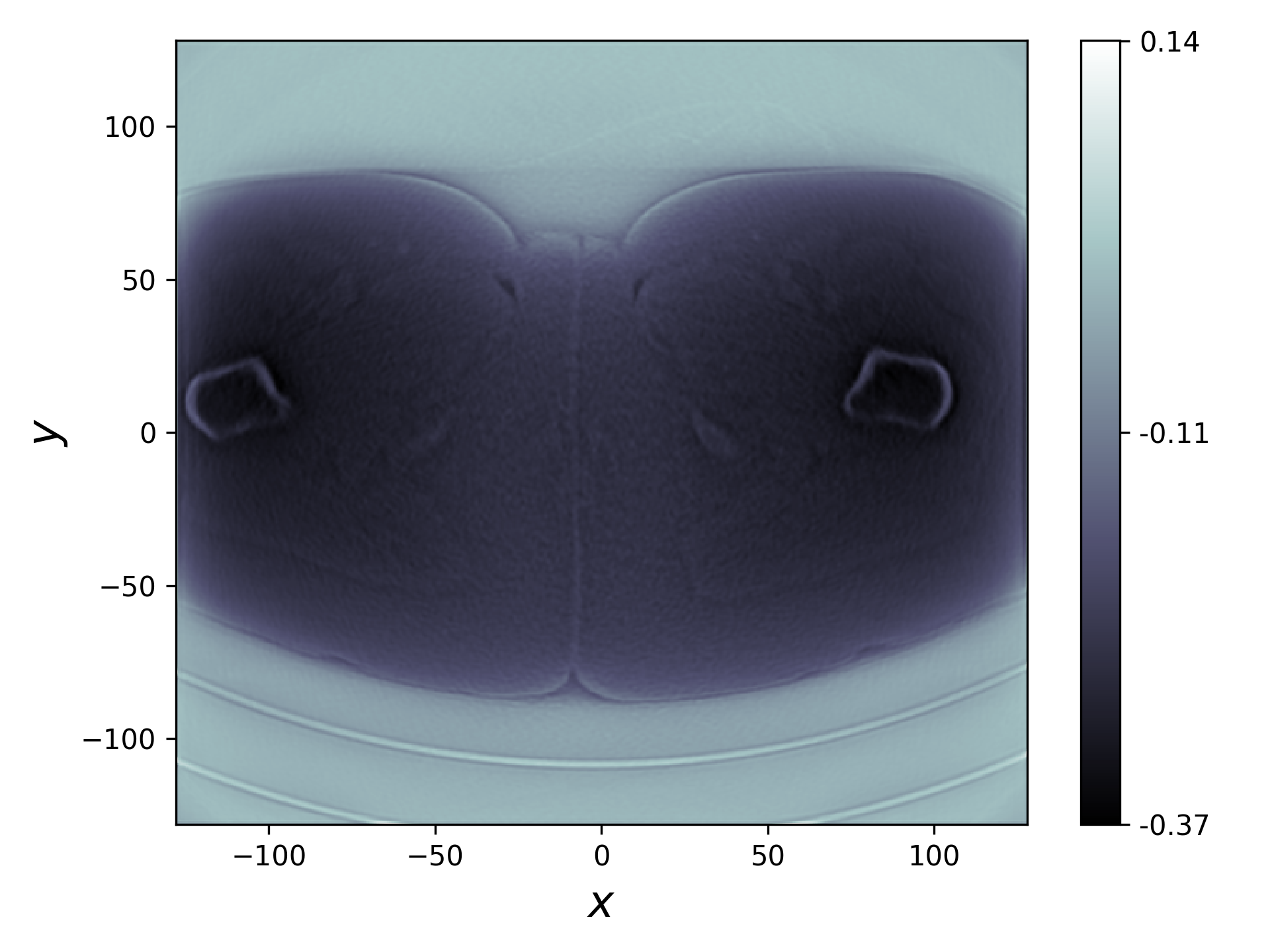}
		\end{subfigure}
		\begin{subfigure}[t]{.32\linewidth}
			\includegraphics[width=\textwidth, trim={23mm 17mm 32mm 6mm}, clip]{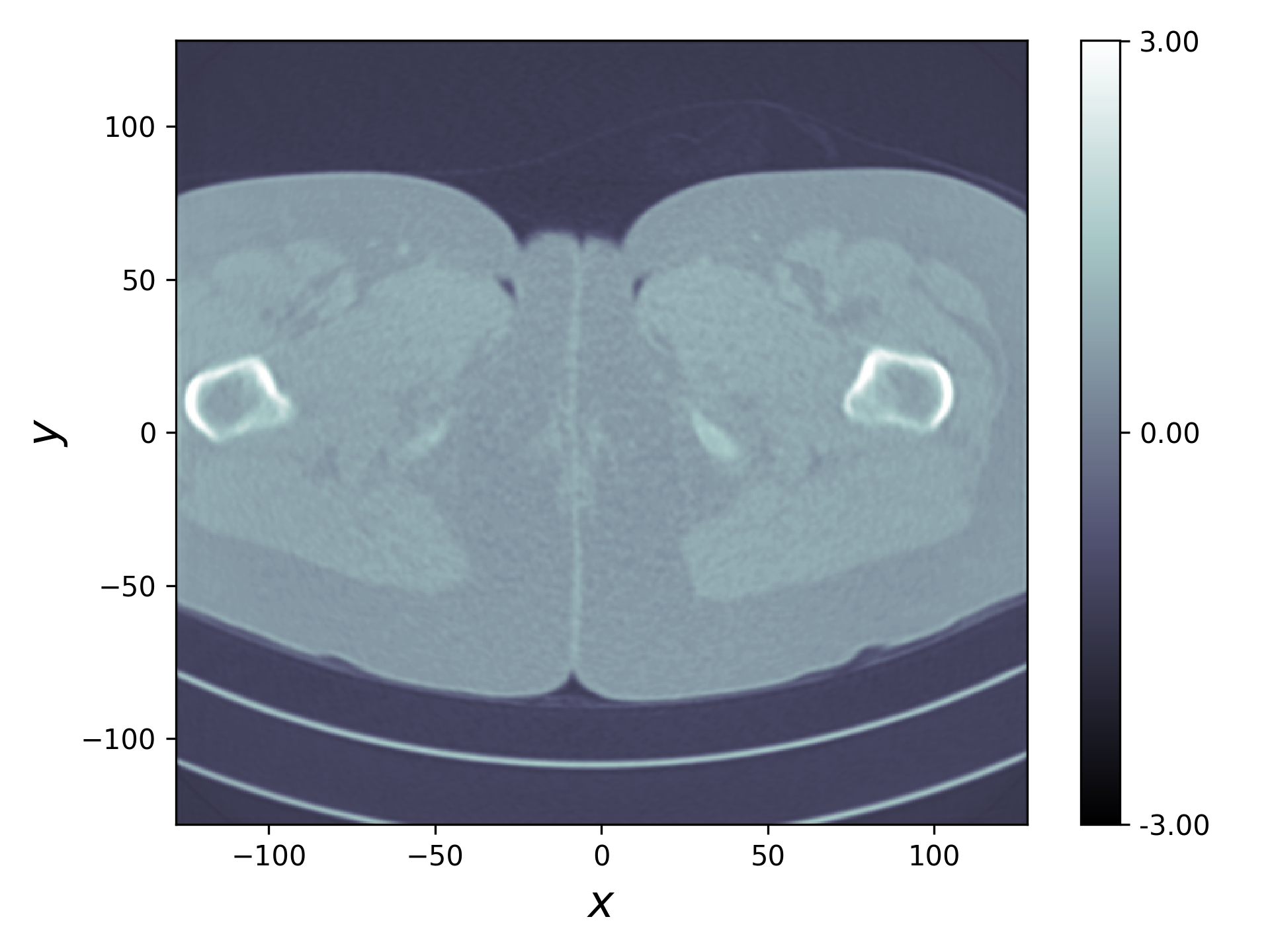}
		\end{subfigure}
		\begin{subfigure}[t]{.32\linewidth}
			\includegraphics[width=\textwidth, trim={23mm 17mm 32mm 6mm}, clip]{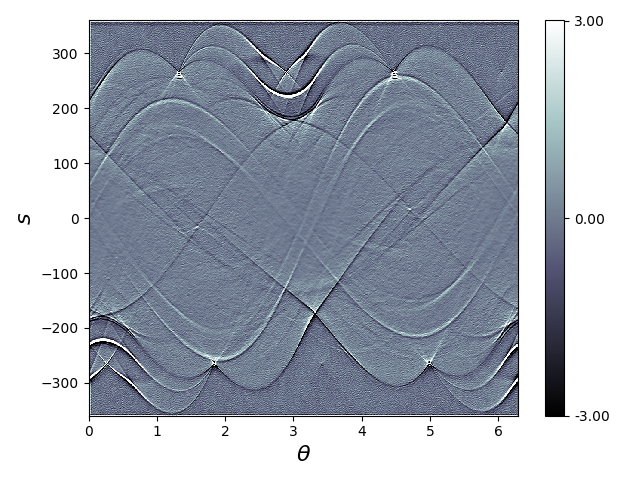}
		\end{subfigure}
		\\ \vspace{1mm}
		\begin{subfigure}[t]{.32\linewidth}	
			\includegraphics[width=\linewidth, trim={23mm 17mm 32mm 6mm}, clip]{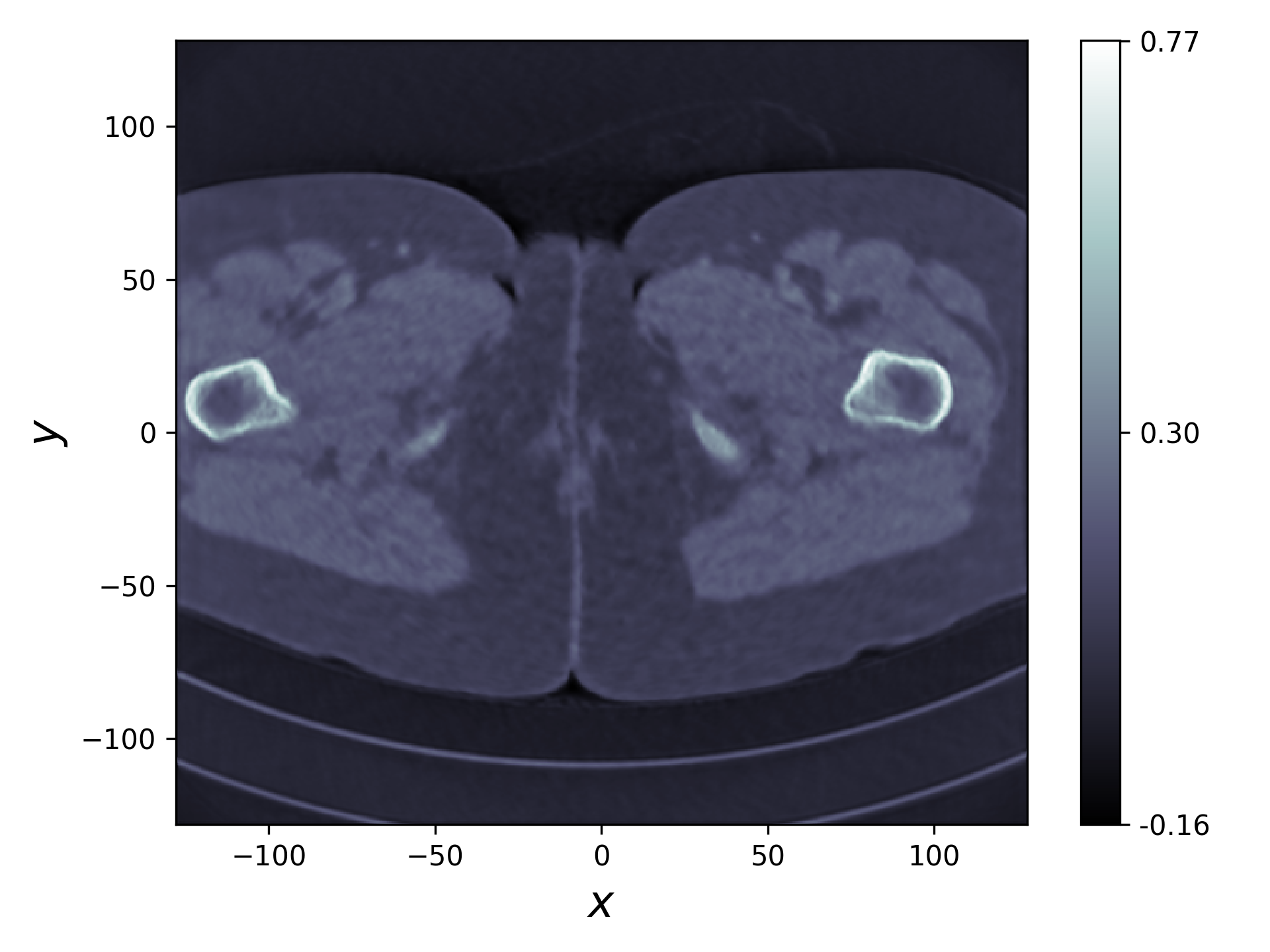}
		\end{subfigure}
		\begin{subfigure}[t]{.32\linewidth}
			\includegraphics[width=\textwidth, trim={23mm 17mm 32mm 6mm}, clip]{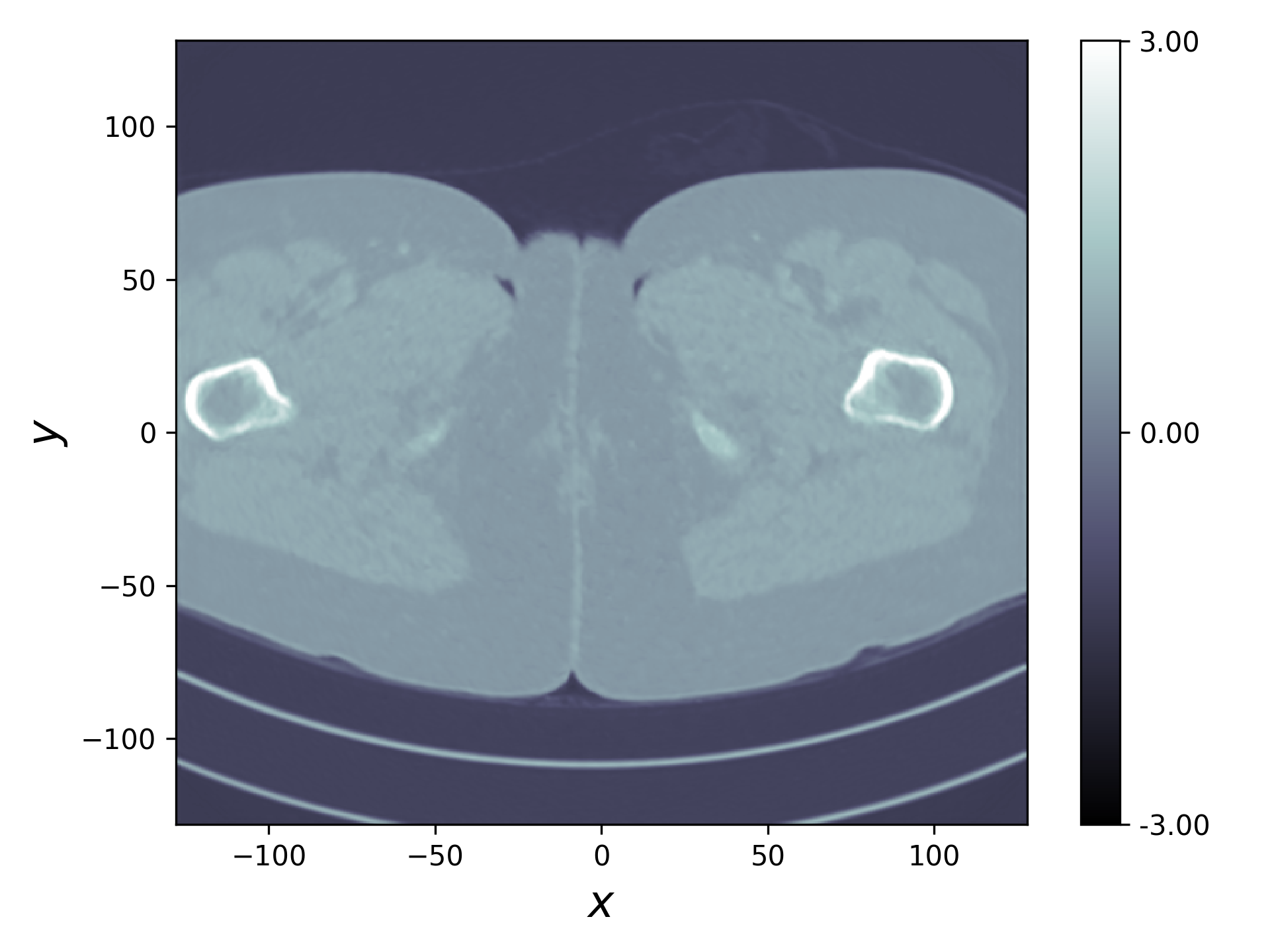}
		\end{subfigure}
		\begin{subfigure}[t]{.32\linewidth}
			\includegraphics[width=\textwidth, trim={23mm 17mm 32mm 6mm}, clip]{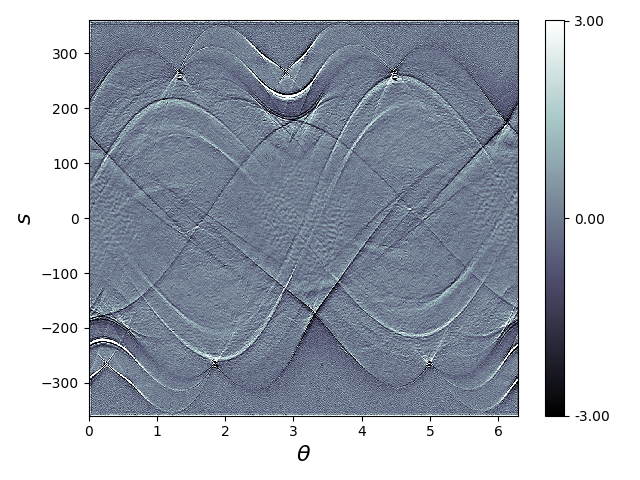}
		\end{subfigure}
		\\ \vspace{1mm}
		\begin{subfigure}[t]{.32\linewidth}	
			\includegraphics[width=\linewidth, trim={23mm 17mm 32mm 6mm}, clip]{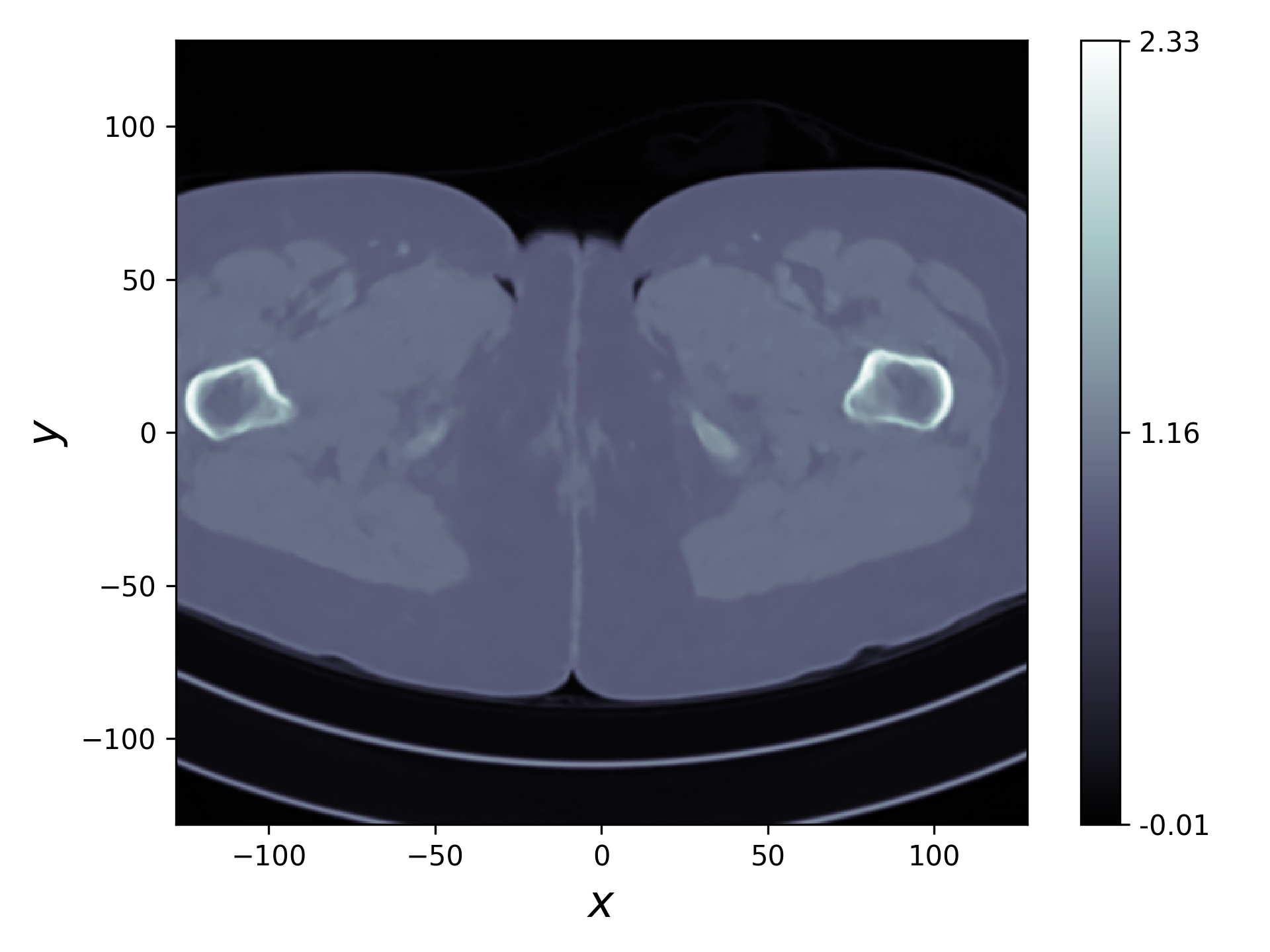}
		\end{subfigure}
		\begin{subfigure}[t]{.32\linewidth}
			\includegraphics[width=\textwidth, trim={23mm 17mm 32mm 6mm}, clip]{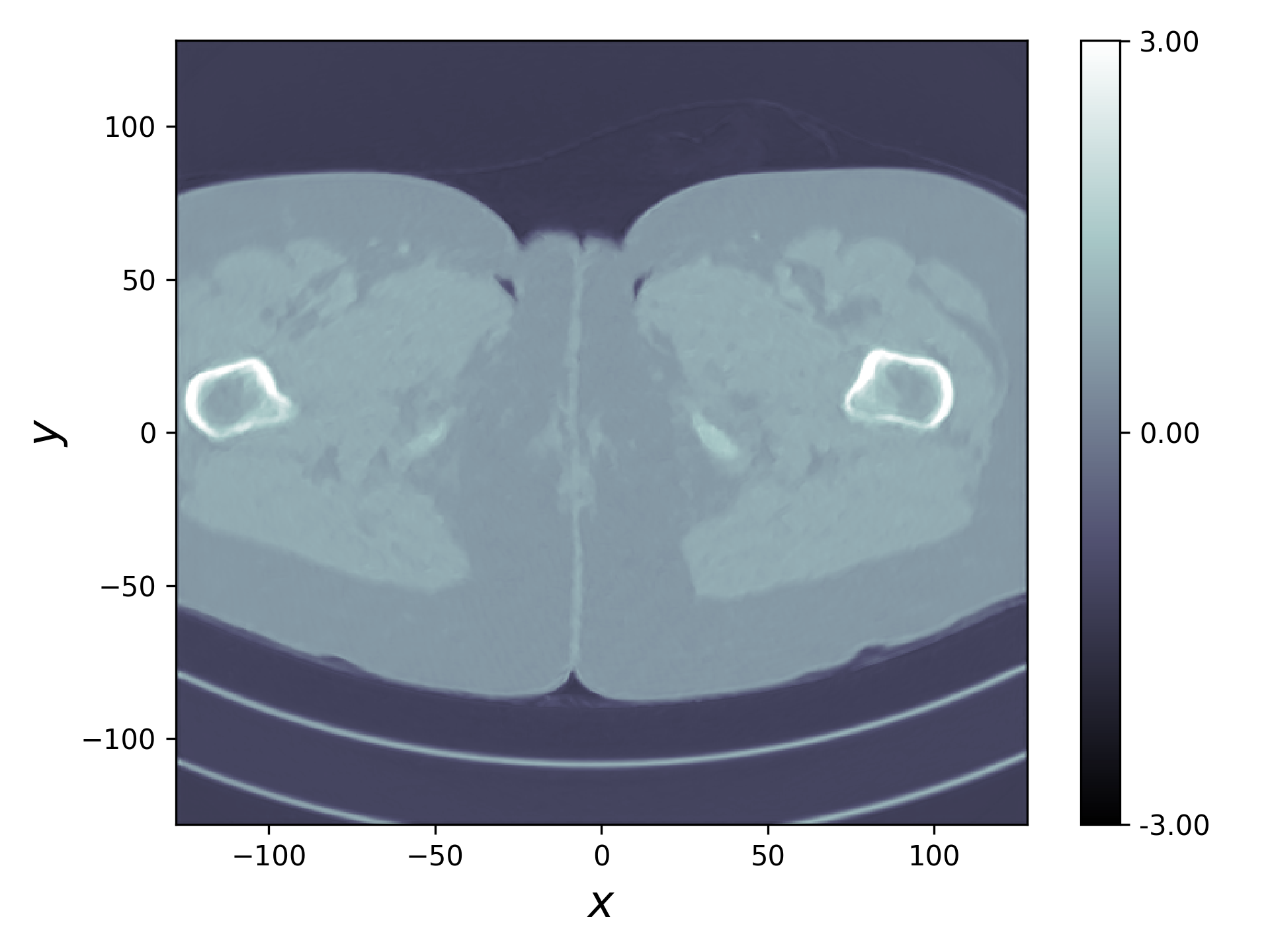}
		\end{subfigure}
		\begin{subfigure}[t]{.32\linewidth}
			\includegraphics[width=\textwidth, trim={23mm 17mm 32mm 6mm}, clip]{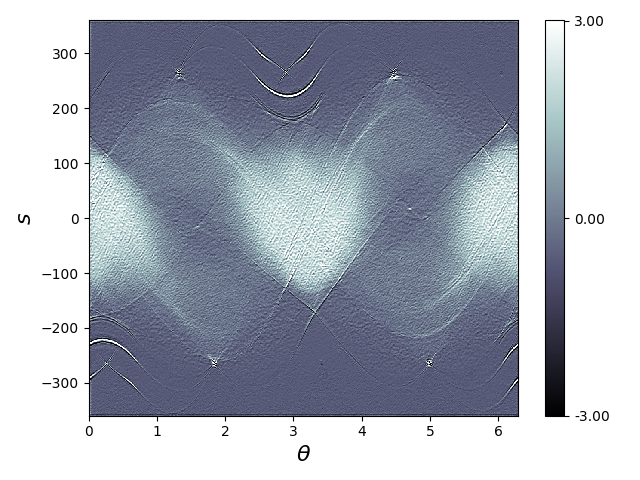}
		\end{subfigure}
		\caption{Iterates 2, 4, 6, 8 and 10 in the Learned Primal-Dual algorithm when reconstructing the human phantoms using a non-linear forward model. Left: Reconstruction ($\primal_i^{(1)}$). Middle: Point of evaluation for the forward operator ($\primal_i^{(2)}$). Right: Point of evaluation for the adjoint of the derivative ($\dual_i^{(1)}$). Windows selected to cover most of the range of the values.}
		\label{fig:iterates}
	\end{figure}

	Since the algorithm is iterative, we can visualize the iterates to gain insight into how it works. In \cref{fig:iterates} we show some iterates with the non-linear forward operator. We note that the reconstruction stays very bad until the 8:th iterate when most structures seem to come in place, but the image is still noisy. Between the 8:th and 10:th iterate, we see that the algorithm seems to perform an edge-enhancing step. It thus seems like the learned iterative scheme works in two steps, first finding the large scale structures and then fine-tuning the details.
	
	Similarly to the edge-enhancement that seems to be performed in the primal space, we note that in the dual space the sinogram that is back-projected seems to be band-pass filtered to exclude both very low and very high frequencies. 	
	
	We note that in the very noisy and under-sampled data used for the ellipse phantoms, the learned algorithms that make use of the forward operator, such as the Learned Gradient, Primal and Primal-Dual algorithms outperform even state of the art post-processing methods by large margins and that in this regimen, \ac{TV} regularization performs relatively well when compared to post-processing methods. This improvement in reconstruction quality when incorporating the forward operator, while still substantial, is not as large for the human phantom in which the data was less noisy.
	
	To explain this, we conjecture that in the limit of highly noisy data where the initial reconstruction as given by e.g. \ac{FBP} becomes very bad, learned schemes that incorporate the forward model and work directly from data, such as the Learned Primal-Dual algorithm, has a considerable advantage over post-processing methods and that this advantage increases with decreasing data quality.
	
	Further along these lines,  note that for the human data the post-processing gives a large improvement in \ac{PSNR} when compared to \ac{TV} regularization, which is not necessarily reflected in the  \ac{SSIM}. On the other hand, the Learned Primal-Dual algorithm improves upon both \ac{PSNR} and \ac{SSIM}. This can be by explained by the learned post-processing being limited by the information content of the \ac{FBP} while the Learned Primal-Dual algorithm works directly with data and is thus limited by the information content of the data, which is greater or equal to that of the \ac{FBP}. In theory, the Learned Primal-Dual algorithm can thus find structures that are not present in the \ac{FBP}, something post-processing methods cannot.
	
	In these experiments we found that while the algorithm seems to handle non-linear forward models well, we did not observe any notable performance improvement by doing so. This may indicate that performing reconstructions on post-log data is preferable.
	
	The structure of the neural network  was not extensively fine-tuned and we suspect that better results could be obtained by a better choice of network  for the learned proximal operators.
	We also observed that the choice of optimizer and learning rate decay had a large impact on the results, and we suspect that further research into how to correctly train learned reconstruction operators will prove fruitful.
	
	Finally, we observe that the reconstructions, while outperforming all of the compared methods with respect to \ac{PSNR} and \ac{SSIM}, suffers from a perceived over-smoothing when inspected visually. We suspect that the particular choice of objective function used in this article, the squared norm \cref{eq:LossStandard}, is a main cause of this and invite future researchers to implement learned reconstruction operators that use more advanced loss functions such as perceptual losses \cite{PerceptualLoss}.
	
	\section{Conclusions}
	
	We have proposed a new algorithm in the family of deep learning based iterative reconstruction schemes. The algorithm is inspired by the \ac{PDHG} algorithm, where we replace the proximal operators by learned operators. In contrast to several recently proposed algorithms, the new algorithm works directly from tomographic data and does not depend on any initial reconstruction.
	
	We showed that the algorithm gives state of the art results on a computed tomography problem for both analytical and human phantoms. For analytical phantoms, it improves upon both classical algorithms such as \ac{FBP} and \ac{TV}, and post-processing based algorithms by at least \unit{6}{\decibel} while also improving the \ac{SSIM}. The improvements for the human phantom were more modest, but the algorithm still improves upon a \ac{TV} regularized reconstruction by \unit{6.6}{\decibel} and gives an improvement of \unit{2.2}{\decibel} when compared to a learned post-processing.
	
	We hope that this algorithm will inspire further research in Learned Primal-Dual schemes and that the method will be applied to other imaging modalities.
	
	\section{Acknowledgements}
	The work was supported by the Swedish Foundation of Strategic Research grant AM13-0049 and Industrial PhD grant ID14-0055. The work was also supported by Elekta.
	
	The authors also thank Dr. Cynthia McCollough, the Mayo Clinic, and the American Association of Physicists in Medicine, and acknowledge funding from grants EB017095 and EB017185 from the National Institute of Biomedical Imaging and Bioengineering, for providing the data necessary for performing comparison using a human phantom.
	
	\newpage
	\bibliographystyle{IEEEtran}
	\bibliography{learning_refs}
\end{document}